\documentclass[11pt]{amsart}
\usepackage{amsmath,amssymb,amsfonts,url}

\newcommand{\Ric}{\operatorname{Ric}}
\renewcommand{\div}{\operatorname{div}}
\newcommand{\Ker}{\operatorname{Ker}}
\newcommand{\Vol}{\operatorname{Vol}}
\newcommand{\D}{\partial}
\newcommand{\gt}{\tilde g}

\newcommand{\norm}[1]{\left| #1 \right|}
\newcommand{\eps}{\varepsilon}
\newcommand{\gh}{\hat g}
\newcommand{\wh}{\hat w}
\newcommand{\bh}{\hat b}

\newcommand{\uh}{\hat u}
\newcommand{\gb}{\breve g}
\newcommand{\bb}{\breve b}
\newcommand{\hb}{\breve h}
\newcommand{\fb}{\breve f}
\newcommand{\rb}{\bar r}
\newcommand{\gbar}{\bar g}
\newcommand{\Rb}{\bar R}
\newcommand{\Rh}{\hat R}
\newcommand{\Gh}{\hat G}
\newcommand{\phih}{{\hat\phi}}
\newcommand{\betah}{\hat\beta}
\newcommand{\p}{\varphi}

\newcommand{\compose}{\operatorname{\scriptstyle\circ}}
\newcommand{\R}{\mathbb R}

\numberwithin{equation}{section}

\newtheorem{thm}{Theorem}
\renewcommand{\thethm}{\kern -.8ex}
\newtheorem{lem}{Lemma}[section]
\newtheorem{rem}{Remark}[section]

\begin{document}
\title[Asymptotic expansion for bubbling solutions]{The profile of bubbling solutions\\
of a class of fourth order\\ geometric equations on 4-manifolds}

\author{Gilbert Weinstein}
\address{Department of Mathematics\\
        University of Alabama at Birmingham\\
        Birmingham, Alabama 35205}
\email{weinstein@uab.edu}

\author{Lei Zhang}
\address{Department of Mathematics\\
        University of Alabama at Birmingham\\
        Birmingham, Alabama 35205}
\email{leizhang@math.uab.edu}
\thanks{Lei Zhang is supported by National Science Foundation Grant 0600275 (0810902)}

\date{\today}

\begin{abstract}
We study a class of fourth order geometric equations defined on a $4$-dimensional compact Riemannian manifold which includes the $Q$-curvature equation.  We obtain sharp estimates on the difference near the blow-up points between a bubbling sequence of solutions and the standard bubble.
\end{abstract}

\maketitle



\section{Introduction}
Let $(M,g)$ be a compact Riemannian manifold.  The \emph{conformal class}
of $g$ consists of all metrics $\gt = e^{2u}g$ for any smooth function $u$.  A central theme in conformal geometry is the study of properties that are common to all metrics in the same conformal class, and the understanding and classification of all the conformal classes.  For this purpose it is often useful to be able to single out a unique representative in each conformal class by imposing some geometric condition.  This usually leads to a conformally covariant geometric equation for the conformal factor $e^{2u}$.  Such equations have attracted much interest in the literature in the past half-century.

In dimension 2, the natural condition to impose is constant Gauss curvature.  The Poincar\'e Uniformization Theorem states that this is always possible: \emph{every compact Riemannian surface is conformal to one with constant Gauss curvature}.  The conformally covariant operator in this case is the Laplacian-Beltrami operator, given in local coordinates by:
\[
    \Delta_g = \frac1{\sqrt{\det g}} \, \D_i\left( \sqrt{\det g}\, g^{ij} \D_j\right),
\]
and the equation for constant curvature is:
\[
    - \Delta_g u + \kappa_g = \kappa_{\gt} e^{2u}
\]
where $\kappa_g$ is the Gauss curvature of $g$, and $\kappa_{\gt}$ is constant.  Using the Gauss-Bonnet theorem, we see that the sign of $\kappa_{\gt}$ is determined by $\chi_M$ the Euler characteristic of $M$:
\[
    2\pi\chi_M = \int_M \kappa_g\, dA_g = \kappa_{\gt} \int_M e^{2u}\, dA_g,
\]
where $dA_g$ is the area element of $g$.  Although this result was originally proved by Poincar\'e using non-PDE methods, there is now a PDE proof, see~\cite{christodoulou} and~\cite{weinstein}.  Furthermore, the operator $\Delta_g$ is conformally convariant $\Delta_{\gt} = e^{-2u} \Delta_g$.

For compact Riemannian manifolds of dimension $n\geq3$, a natural generalization is to impose constant scalar curvature.  This leads to the Yamabe problem: \emph{given a compact Riemannian manifold $(M,g)$ of dimension at least $3$, find a metric conformal to $g$ with constant scalar curvature}.  This was also eventually answered in the affirmative, see~\cite{yamabe,trudinger,aubiny,schoen}.  The corresponding operator is now the conformal Laplacian $L_g=\Delta_g-c_nR_g$, $c_n=(n-2)/4(n-1)$, and the equation for constant scalar curvature is the Yamabe equation:
\[
    L_g \phi = \epsilon \phi^{(n+2)/(n-2)}
\]
where $\phi^{4/(n-2)}=e^{2u}$, and $\epsilon$ is $1$, $0$, or $-1$.  The operator $L_g$ also has a conformal covariant property:
\[
    L_{\gt} \phi = \phi^{-(n+2)/(n-2)} L_g \phi,
\]
where $\gt = \phi^{4/(n-2)}g$.

In $4$-d, another problem analogous to the $2$-d case arises from imposing the condition of constant $Q$-curvature:
\[
    Q_g = -\frac 1{12}(\Delta_g R_g - R_g^2 + 3|\Ric_g|^2).
\]
The natural question is the same: \emph{given a $4$-d compact Riemannian manifold $(M,g)$, is there a metric $\gt = e^{2u} g$ in the conformal class of $g$ with constant $Q$-curvature?}  The $Q$-curvature of the metric $\gt$ is given by:
\[
    P_g u + 2Q_g = 2Q_{\tilde g}e^{4u},
\]
where $P_g$ is the Paneitz operator:
\[
    P_g u = \Delta_g^2 u + \div_g\left( \left(\frac23 R_gg-2\Ric_g\right) \nabla u
    \right).
\]
Integrating with respect to the volume element $dV_g$, it is easy to see that the quantity:
\[
    k_P = \int_M Q_g\, dV_g
\]
is a conformal invariant, i.e., it is constant in the conformal class of $g$.  Furthermore, we also have a Gauss-Bonnet formula:
\[
    \int_M \left( Q_g + \frac18 \norm{W_g}^2 \right) dV_g = 4\pi^2 \chi_M,
\]
where $W_g$ is the Weyl tensor of $g$ given in local coordinates as:
\[
    W_{ijkl} = R_{ijkl} - \frac2{n-2}(g_{i[k} R_{l]j} - g_{j[k} R_{l]i})
        + \frac{2}{(n-1)(n-2)}R g_{a[k}g_{l]j}.
\]
We note that $W_g$ is pointwise conformally invariant $W_{\gt}=W_g$, and the operator $P_g$ is conformally covariant:
\begin{equation}    \label{jan22e2}
    P_{\gt} f = e^{-4u} P_g f.
\end{equation}

We use \eqref{P} to denote the assumption:
\[  \label{P}
    \Ker(P_g)=\{\text{constants}\}.
    \tag{$P$}
\]
We remark that \eqref{P} is often satisfied. For example, Gursky
in~\cite{gursky} proved that if $(M,g)$ has
non-negative Yamabe invariant $Y_g\geq0$, and satisfies $k_P\ge
0$, then \eqref{P} holds and $P_g\ge 0$; see
also~\cite{gurskyviaclovsky} where the assumptions are weakened to
$Y_g\geq0$ and $k_P+Y_g^2/6>0$.

Chang and Yang proved in~\cite{changyang1} that if $k_P < 8\pi^2$, $P_g\ge 0$ and \eqref{P} holds, then there is a conformal metric $\gt$ whose $Q$-curvature is constant. In~\cite{djadli2}, Djadli and Malchiodi extended this existence result assuming only that \eqref{P} holds and $k_P\neq 8\pi^2N$ for any positive integer $N$. An essential ingredient in this existence result is an a priori bound: if $k_P \neq 8\pi^2N$ for any positive integer $N$, then any sequence of solutions of the prescribed $Q$-curvature equation is uniformly bounded.  In fact, this a priori estimate can be extended to the following more general equation in the same class:
\begin{equation}
 \label{jan22e1}
    P_g u + 2b = 2h e^{4u},
\end{equation}
where $b$ is a smooth function.  Note that if $b=Q_g$, then $h$ is the $Q$-curvature of the conformal metric $e^{2u}g$.  Assuming $h_k\to h_0$, $h_k\geq c_0>0$, and $b_k\to\ b_0$, Druet and Robert in~\cite{druet1} showed that
any sequence of solutions $\{u_k\}$ of~\eqref{jan22e1} with $h=h_k$ and $b=b_k$ is uniformly bounded, provided $\int_Mb_0\neq 8\pi^2N$, see also Malchiodi~\cite{mal1}, .

However, bubbling can occur when $\int_M b_0\, dV_g = 8\pi^2 N$ for some positive integer $N$. A precise understanding of this bubbling phenomenon is required if progress is to be made on the existence problem.  The study of the blow-up profile and other blow-up phenomena for the Paneitz operator and other $4$-th order elliptic equations has attracted much interest recently; see for example~\cite{adimurthi,brendle,changqingyang,djadli1,fefferman,gurskyviaclovsky2,lililiu,malchiodi,malchiodi2,ndiaye,qing, struwe,wei,wei2}.

Let $\{u_k\}$ be a sequence of solutions of~\eqref{jan22e1} with $h=h_k$, and $b=b_k$.  We say that this is a \emph{bubbling sequence} if $\sup |u_k| \to\infty$.  In~\cite{druet1}, Druet and Robert studied bubbling sequences of solutions of~\eqref{jan22e1} and obtained
some asymptotic estimates on the behavior near the blow-up points.
We will throughout make the following assumptions on the coefficients $b_k$ and $h_k$:
\[  \label{b}
    \| b_k-b_0 \|_{C^1(M)} \to 0, \qquad
     \|h_k - h_0\|_{C^2(M)} \to 0, \qquad h_k\geq c_0.
    \tag{$b,h$}
\]
It follows immediately that $\|b_k\|_{C^1(M)}\leq C_0$, and $\|h_k\|_{C^2(M)}\leq C_0$ for some constant $C_0$ independent of $k$.
We let $G$ denote the Green's function for the Paneitz operator:
\begin{equation}    \label{green}
    f(\xi)-\bar f_g=\int_M G(\xi,\eta)\, P_gf(\eta)\, dV_g(\eta),  \qquad \int_M G(\xi,\eta)\, dV_g(\eta) = 0,
\end{equation}
where $\bar f_g=\Vol_g(M)^{-1} \int_M f\, dV_g$ is the mean value of $f$.  The asymptotics of this Green's function are studied in the Appendix.  Now, for $k=0,\dots$, let
\begin{equation} \label{mar26e1}
    \phi_k(\xi) = 2 \int_M G(\xi,\eta)\, b_k(\eta)\, dV_g(\eta).
\end{equation}
Since $\{u_k\}$ is a bubbling sequence, it follows immediately that $\int_M b_0\, dV_g = 8 N\pi^2$ for some positive integer $N$.  Druet and Robert proved that passing to a subsequence, there is a finite set $S=\{p_1,..,p_N\}$ such that:
\[
    u_k - \bar u_k \to 16\pi^2 \sum_i^N G(p_i,\cdot) - \phi_0 \qquad
    \text{in $C^4_{loc}(M\setminus S)$},
\]
Let $\beta$ be the regular part of the Green's function:
\begin{equation}    \label{beta}
    G(\xi,\eta)= - \frac{1}{8\pi^2}\, \chi(r)\,
    \log d_g(\xi,\eta) + \beta(\xi,\eta).
\end{equation}
Here $\chi$ is a cut off function supported in a small
neighborhood of $\xi$, and $r=d_g(\xi,\eta)$. They also proved
that for $i=1,\dots,N$:
\[
    64\pi^2 \nabla_{2} \beta(p_i,p_i) + 64\pi^2 \sum_{j\neq i}\nabla_{1} G(p_i,p_j)
    - 4\nabla \phi_0(p_i) = - \frac{\nabla h(p_i)}{h(p_i)},
\]
where $h$ is the limit of $h_k$ as $k\to \infty$, and $\nabla_1$, $\nabla_2$ denote the derivatives with respect to the first and second variables respectively.  In this article we will continue this line of investigation and derive more precise asymptotic estimates for the behavior of such solutions.  We define the \emph{standard bubble} at $p$:
\[
    U_{p,\eps,H}(\xi) = - \log\left( \eps + \frac{\sqrt{H}\, d_g(p,\xi)^2}{4\sqrt3\,\eps} \right).
\]
We will also adopt the following notation.  For $k$ large enough, there are $N$ points $\{q_{ik}\}$ such that $q_{ik}\to p_i$ and $u_k(q_{ik})\to\infty$.  Let $H_{ik}=h_{k}(q_{ik})$, $\eps_{ik}=e^{-u_k(q_{ik})}$, and $U_{ik} = U_{q_{ik},\eps_{ik},H_{ik}}$.

\begin{thm}
Let $\{u_k\}$ be a bubbling sequence of solutions on $M$.  Then passing to a subsequence, there is a constant $\delta>0$ such that for any fixed $\tau\in (0,1)$, there exists a constant $C_1 = C_1(N,g,c_0,C_0,\tau)$ such that:
\begin{equation}    \label{MainEst}
    \left| u_k(\xi) - U_{ik}(\xi) \right|
    \le C_1 d_{g}(q_{ik},\xi)^{\tau},
\end{equation}
in $B(q_{ik},\delta)$, and such that for $i=1,\dots, N$ we have:
\begin{multline}    \label{vrate}
    \left| 64\pi^2\nabla_2\beta(q_{ik},q_{ik}) + 64\pi^2\sum_{j\neq i}\nabla_1 G(q_{ik},q_{jk})
    \right. \\ \left.
    - 4\nabla \phi_k(q_{ik})+\frac{\nabla h_k(q_{ik})}{h_k(q_{ik})}\right|
    \le C_1 e^{-\tau u_k(q_{ik})/2}.
\end{multline}
\end{thm}

Our approach is motivated by Lin and Wei's work~\cite{linwei1}
from which one can easily derive an $O(1)$ bound
in~\eqref{MainEst} (i.e. $|u_k-U_{ik}|\le
C$ near $p_i$) provided $(M,g)$ is locally conformally flat; see also
~\cite{xu} which uses a completely different approach. Our result
removes the hypothesis of local conformal flatness and also
improves the estimate near the blow-up points. We hope our
approach can be fine-tuned to yield better yet estimates as
required to handle the existence question posed above.

A major difficulty when trying to prove a priori estimates for solutions of fourth order elliptic equations is the lack of a maximum principle. In order to remedy this, Lin and Wei devised a strategy based on the Pohozaev identity.  We adapt this approach to the case in which the manifold is not necessarily locally conformally flat, making use of \emph{conformal normal coordinates}.  These are normal coordinates for a metric $\gh$ in the conformal class of $g$ for which $\det(\gh)=1$.  The existence of such a metric is proved in~\cite{cao}.  Although, we used this result for the sake of simplicity, our proof only relies on the weaker concept already introduced in the solution of the Yamabe problem where one only requires $\det(\gh)=1$ to hold to high enough order in the distance from the center of the ball under consideration, see~\cite{leeparker}.

We now briefly sketch the outline of the paper and the proof of our Theorem.

In Section~\ref{sec:02}, we prove the $O(1)$ estimate.  We use the Green's representation formula, together with rough estimates from~\cite{druet1}, to write long range asymptotic formulas for the rescaled solution $v_k$ and its derivatives in terms of the concentration of energy $\alpha_k$ near the singular point, i.e. within a carefully chosen radius $l_k = -\eps_k\log\eps_k$, where $\eps_k$ is related to the maximum of $u_k$, see~\eqref{alpha}.  These are then substituted into an asymptotic Pohozaev identity, and after estimating the higher order terms, we obtain an asymptotic formula for the energy $\alpha_k \approx 16\pi^2$.  When substituted back into the asymptotic formula for $v_k$, this yields a long range $O(1)$ bound.  Finally, we use standard estimates in the interior, and then these long range and interior estimates in conjunction with the maximum principle in the mid-range.

In Section~\ref{sec:03}, we prove~\eqref{MainEst} by contradiction.  We divide the argument into two cases, depending upon whether an appropriately weighted supremum runs off to infinity or remains in a bounded region along a subsequence.  In the first case, we use the Green's representation formula and a comparison between the geometric and Euclidean distances to reach a contradiction.  In the second case, we show that the difference between the appropriately rescaled solutions and the standard bubble converges to a solution of the linearized equation which we can then show vanishes thanks to a lemma of Lin and Wei from~\cite{linwei1}, again leading to a contradiction.

Finally, in Section~\ref{sec:04}, we use our estimate~\eqref{MainEst} in the Pohozaev Identity over a ball of radius $\eps_k^{-1/2}$ to obtain a Euclidean version of the vanishing rate.  We then translate this result into the original metric $g$ and prove~\eqref{vrate}.

The Appendix deals with delicate estimates for the Green's function, an asymptotic comparison between the geodesic distance and the Euclidean distance in conformal normal coordinates, some well known curvature and metric derivatives computations in conformal normal coordinates, and a proof of the asymptotic Pohozaev identity.

\section{The $O(1)$ estimate}
\label{sec:02}

In this section we derive the $O(1)$ estimate, i.e., we show that
\[
    |u_k(\xi)-U_{ik}(\xi)|\le C, \qquad \text{for $\xi\in B(q_{ik},\delta)$.}
\]
This estimate has been established by Lin-Wei \cite {linwei1} for
locally conformally flat manifolds; see also ~\cite{xu} for a
completely different proof. Furthermore, we remove the assumption
that $(M,g)$ is locally conformally flat.

 Our first step
is to rescale the solutions, and use the Green's representation
formula~\eqref{green} to derive the long range asymptotic
formulas~\eqref{nov2e1}--\eqref{nov2e6}.

In~\cite{druet1}, Druet and Robert prove that the singular set $S$ consists of only finitely many points $\{p_1,\dots,p_N\}$ and these are separated uniformly in $k$ by a positive distance.  Without loss of generality we will focus in this section on $p_1$, and to simplify the notation, we will omit the subscript $1$, so that we now consider a sequence of points $q_k\in M$ where $u_k$ has a local maximum $u_k(q_k)\to\infty$ and $q_k\to p$ as $k\to\infty$.

According to~\cite{cao}, we can find function $\wh_k$ defined on
$M$, such that in a neighborhood $B(q_k,\delta_1)$ of $q_k$,
$\delta_1>0$, we have $\det(\gh_k) = 1$ in the normal coordinates
of the conformal metric $\gh_k=e^{2\wh_k}g$.  We refer to these
coordinates as \emph{conformal normal coordinates}. We point out
that $\det(\gh_k)\approx1$ to high enough order would be
sufficient for our purpose, but we use Cao's result since it
simplifies the proof.  We also choose $\delta_1$ small enough so
that $\delta_1<\operatorname{inj}(M)/10$ and $\delta_1<d/10$ where
$d$ is the minimum distance between any two points in the singular
set $S$. Using the conformal covariance property of
$P_g$~(\ref{jan22e2}), we obtain that the function $\uh_k =
u_k-\wh_k$ satisfies\label{bkhat}
\[
    P_{\gh_k}\uh_k+2\bh_k=2h_ke^{4\uh_k}.
\]
where $2\hat b_k=P_{\gh_k}\wh_k+2b_ke^{-4\wh_k}$.  We remark that if $b_k=Q_g$ then $\bh_k=Q_{\gh_k}$.  We also note that $\wh_k(\xi)=O(d_g(\xi,q_k)^2)$ in a neighborhood of $q_k$, hence all the terms coming from $\wh_k$ can be absorbed on the right-hand side of~\eqref{MainEst}.  We have the following estimates, also proved in~\cite{druet1}:
\begin{equation} \label{oct31e2}
\begin{gathered}
    \eps_k^{(1-\nu)}d_{\gh_k}(\xi,q_k)^{\nu}e^{\uh_k(\xi)} \le C_{\nu}, \qquad 1\le \nu<2 \\
     \norm{D^j\uh_k(\xi)} \le C(d_{\gh_k}(\xi,q_k))^{-j},\qquad j=1,2,3,
\end{gathered}
\end{equation}
where $\norm{D^j\uh_k(\xi)} = \sum_J \norm{D^J \uh_k(\xi)}$ and the sum is over all multi-indices $J$ of order $j$, and $\eps_k=e^{-\uh_k(q_k)}$.  We now rescale the solutions $\uh_k$, using a blow-up of the neighborhood of the point $q_k$.  Define the map $\p_k\colon  B(0,\delta_1\eps_k^{-1})\to B(q_k,\delta_1)$ by $\p_k\colon y \mapsto \eps_k y$, where on the right-hand side we are using conformal normal coordinates on $B(q_k,\delta_1)$.  We use the notation $\fb = \p_*f = f\compose\p$ to denote the pull-back of a function $f$ defined on $B(q_k,\delta_1)$, and we let $\gb_k = \eps_k^{-2}\p_* g_k$ be the blow-up metric, i.e., a rescaling of the pull-back metric.  We define:
\[
    v_k = \breve{\uh}_k + \log\eps_k,
\]
and note that $v_k(0)=0$.  It follows from~\eqref{jan22e2} that $v_k$ satisfies:
\begin{equation}    \label{oct31e1}
    P_{\gb_k} v_k + 2\eps_k^4 \bb_k = 2 \hb_k e^{4v_k},
    \qquad  \text{in $B(0,\delta_1\eps_k^{-1})$.}
\end{equation}
The estimates~\eqref{oct31e2} now read:
\begin{gather}
    \label{oct31e3} |v_k(y)|\le (-2+\mu)\, \log(1+|y|)+C(\mu), \quad |y|\le \delta_1\eps_k^{-1} \\
    \label{nov29e8} |D^jv_k(y)|\le C(1+|y|)^{-j}, \quad j=1,2,3.
\end{gather}
where $\mu\in (0,1)$.

Let $\Gh_k$ be the Green's function for $P_{\gh_k}$.  Then, we have:
\[
    \uh_k(\xi) = \overline{\uh_k} + 2 \int_M \Gh_k(\xi,\eta) h_k(\eta)e^{4\uh_k(\eta)}\, dV_{\gh_k}(\eta)
    - 2 \int_M \Gh_k(\xi,\eta)\bh_k(\eta)dV_{\gh_k}(\eta).
\]
where $\overline{\uh_k}$ is the mean value of $\uh_k$. Decompose
$\Gh_k$ into a principal part and a regular part as follows:
\[
    \Gh_k(\xi,\eta) = - \frac1{8\pi^2} \, \chi(r)\, \log d_{\gh_k}(\xi,\eta)
    + \betah(\xi,\eta) = H(\xi,\eta) + \betah(\xi,\eta)
\]
where $\chi=1$ on $B(q_k,\delta_1)$, $\chi=0$ on $M\setminus B(q_k,2\delta_1)$.  We have
\begin{equation}    \label{uhat}
    \uh_k(\xi) = \overline{\uh_k} + 2 \int_{M} H(\xi,\eta)\, h_k(\eta)\, e^{4\uh_k(\eta)}\,
    dV_{\gh_k}(\eta)+\phih_k(\xi)
\end{equation}
where
\begin{equation}    \label{feb15e1}
    \phih_k(\xi) = 2 \int_M \betah(\xi,\eta)\, h_k(\eta)\, e^{4\uh_k(\eta)}\, dV_{\gh_k}(\eta)
    -2 \int_M \Gh_k(\xi,\eta)\,\bh_k(\eta)\,dV_{\gh_k}(\eta).
\end{equation}
Note that since $\det(\gh_k) = 1$ in $B(q_k,\delta_1)$, we have $dV_{\gh_k}(\eta)=d\eta$ in $B(q_k,\delta_1)$.
Taking the difference of~\eqref{uhat} evaluated at $\xi$ and $p_k$, we get:
\begin{multline}    \label{oct31e4}
    \uh_k(\xi) - \uh_k(p_k) = \frac1{4\pi^2} \int_M
    \log\left(\frac{|\eta-q_k|}{d_{\gh_k}(\xi,\eta)}\right)\, \chi(r)\, h_k(\eta)\, e^{4\uh_k(\eta)}\,
    dV_{\gh_k}(\eta)\\
    + \phih_k(\xi)-\phih_k(q_k).
\end{multline}
Here we have used the fact that since the coordinates are normal $d_{\gh_k}(\eta,q_k) = |\eta-q_k|$.  Thanks to the cut-off function $\chi$, we can now replace the integral over $M$ by an integral over $B(q_k,2\delta_1)$, and after rescaling, we now obtain:
\begin{multline} \label{oct31e5}
    v_k(y) = \frac 1{4\pi^2} \int_{B(0,2\delta_1\eps_k^{-1})}
    \log\left(\frac{|z|}{d_{\gb_k}(z,y)}\right)\, \chi(\eps_k r)\, h_k(\eps_k z)\, e^{4v_k(z)}\, dz \\
    + \phih_k(\eps_k y) - \phih_k(0).
\end{multline}

Let
\[
    l_k = - \eps_k\log\eps_k, \qquad L_k=-\log\eps_k,
\]
and define:
\begin{equation}    \label{alpha}
    \alpha_k = 2\int_{B(q_k,\delta_1)} h_k(\eta)\, e^{4\uh_k(\eta)}\, dV_{\gh_k}
\end{equation}
By (\ref{oct31e3}) and the fact that $\det(\gh_k)=1$, one sees easily that:
\begin{equation}
    \label{jan11e1} \alpha_k = 2\int_{B(q_k,l_k)} h_k(\eta)\, e^{4\uh_k(\eta)}\, d\eta + O(L_k^{-3}).
\end{equation}

As in~\cite{linwei1}, the representation formula~\eqref{oct31e5} implies the following long range asymptotic formulas for $v_k$ and its derivatives:
\begin{align}
    \label{nov2e1}
    v_k(y) &= - \frac{\alpha_k}{8\pi^2}\, \log |y| + O(1),\qquad L_k\le |y|\le \delta_1\eps_k^{-1} \\[1ex]
    \label{nov2e2}
    \partial_r v_k(y) &= -\frac{\alpha_k}{8\pi^2}\, L_k^{-1} + O(L_k^{-2}), \qquad |y|=L_k \\[1ex]
    \label{nov2e3}
    \partial_r\bigl(r\partial_r v_k(y)\bigr) &= O(L_k^{-2}), \qquad |y|=L_k \\[1ex]
    \label{nov2e5}
    \Delta v_k(y) &= -\frac{\alpha_k}{4\pi^2}\, L_k^{-2} + O(L_k^{-3}), \qquad |y|=L_k\\[1ex]
    \label{nov2e6}
    \partial_r \Delta v_k(y) &= \frac{\alpha_k}{2\pi^2}\, L_k^{-3}+ O(L_k^{-4}), \qquad |y|=L_k.
\end{align}
Note that while the asymptotic formula is required in the whole range $L_k\le |y|\le \delta_1\eps_k^{-1}$ for $v_k$, it is only required on $|y|=L_k$ for the derivatives.

Since the proof of these estimates is similar to the one in~\cite{linwei1}, we will only briefly sketch the argument pointing out the main differences, a major one being the difference between the Euclidean and Riemannian distance.  It follows from (\ref{oct31e5}) that for $|y| \ge L_k$:
\begin{equation}    \label{nov30e2}
    v_k(y) =  - \frac1{4\pi^2} \int_{B(0,\delta_1 \eps_k^{-1})} \log d_{\gb_k}(y,z)\,
    h_k(\eps_k z)\, e^{4v_k(z)} \, dz + O(1).
\end{equation}
Moreover for any multi-index $J$ of order $j=1,2,3$, and $|y|=-\log {\eps_k}$, we have:
\begin{equation}    \label{nov30e3}
    D^J v_k(y) =-\frac 1{4\pi^2}\int_{B(0,\delta_1\eps_k^{-1})}
    D^J_{y}(\log d_{\gb_k}(y,z))\, h_k(\eps_k z)\, e^{4v_k(z)}\, dz + O(\eps_k^{j}).
\end{equation}
We now divide the domain of integration in~\eqref{nov30e2} and~\eqref{nov30e3} into three subsets $B(0,\delta_1 \eps_k^{-1}) = \Omega_1 \cup \Omega_2 \cup \Omega_3$, \label{omega}where:
\[
    \Omega_{1} = \{ |z| < |y|/2 \}, \quad
    \Omega_{2} = \{ |z-y|< |y|/2 \},\quad
    \Omega_{3} = B(0,\delta_1\epsilon_k^{-1}) \setminus (\Omega_1 \cup \Omega_2).
\]
We also use the following approximations of the distance $d_{\gb_k}$ and its derivatives by their Euclidean counterparts\footnote{These approximations are proved in the Appendix}:
\begin{gather}
    \label{dgapprox}
    \log d_{\gb_k}(y,z) - \log |y-z| = O(1), \qquad z\in B(0,\eps_k^{-1}\delta_1/2)\\
    \label{mar19e1}
    \left| D^j \bigl( \log d_{\gb_k}(y,z) - \log|y-z| \bigr) \right| \le
    C \eps_k^2|y|^{2-j}, \qquad z\in \Omega_{1}.
\end{gather}
Over $\Omega_2\cup\Omega_3$, the integral~\eqref{nov30e2} can be estimated simply by using~\eqref{oct31e3} and the approximation~\eqref{dgapprox} leading to:
\[
    \int_{\Omega_2\cup \Omega_3}\log d_{\gb_k}(y,z)\, \hb_k(z) \, e^{4v_k(z)}\, dz = O(|y|^{-4+\mu_1}),
    \qquad \text{$\mu_1>0$ small.}
\]
In order to capture the asymptotics of the integral~\eqref{nov30e2} over $\Omega_1$, we again use the approximation~\eqref{dgapprox} to reduce the calculation to the Euclidean case, so that~\eqref{nov2e1} then follows with the help of~\eqref{jan11e1}.

Similarly the estimate of~\eqref{nov30e3} over $\Omega_{2}\cup \Omega_{3}$ can be obtained from the bounds:
\[
    \left| D^j \log d_{\gb_k}(y,z) \right|\le C|y-z|^{-j},\qquad j=1,2,3,\qquad
    z\in \Omega_{2}\cup\Omega_{3},
\]
leading to:
\[
    \int_{\Omega_{2}\cup\Omega_{3}} \left| D^j(\log d_{\gb_k}(y,z)) \right|\,
    \hb_k(z)\, e^{4v_k(z)}\, dz = O(|y|^{-4-j+\mu_1}),
\]
while the estimation of these integrals over $\Omega_{1}$ requires the more precise approximation~\eqref{mar19e1}.

We now use the long range estimates~\eqref{nov2e1}--\eqref{nov2e6} in the following Pohozaev identity\footnote{The proof of this identity can be found in the Appendix.} over the ball $\Omega=B(q_k,l_k)$:
\begin{multline*}
    \int_{\Omega} (2 h e^{4\uh} + \frac12 \xi^i \partial_i h e^{4\uh})
    = \int_{\partial \Omega} \left( \frac12 \xi^i \nu_i h e^{4 \uh}
    - \nu_j \xi^m \gh^{ij} \partial_i (\Delta_{\gh} \uh)  \partial_m \uh
    \right. \\ \left.
    {} + \nu_j \gh^{ij} \Delta_{\gh}\uh \partial_i \uh
    + \nu_j \xi^m \gh^{ij} \Delta_{\gh} \uh \partial_{im} \uh
    - \frac12 \xi^i \nu_i (\Delta_{\gh}\uh)^2 \right) \\
    + \int_{\Omega}\left(\Delta_{\gh} \uh \partial_i \gh^{ij} \partial_j \uh
    + \xi^m \Delta_{\gh} \uh \partial_{im} \gh^{ij} \partial_j \uh
    + \xi^m \Delta_{\gh} \uh \partial_m \gh^{ij} \partial_{ij} \uh
    - 2 \bh \xi^i \partial_i \uh \right) \\
    + 2 \int_{\partial\Omega} \left( \Rh_{ij,l}(0) \xi^l \xi^m \nu_i \partial_j \uh \partial_m \uh
    + O(r^3) |D\uh|^2 \right) \\
    - \int_{\Omega} \left( 2\Rh_{ij,l}(0) (\xi^l \partial_j \uh \partial_i \uh
     + \xi^m \xi^l \partial_j \uh \partial_{im} \uh) + O(r^2) |D\uh|^2 + O(r^4) |D^2\uh| \right).
\end{multline*}
Here, we used the conformal normal coordinates $\xi^i$
on this ball, we denoted $r=|\xi|$, and denoted the unit normal to the boundary by
$\nu_i$.  Furthermore, to simplify the notation, we suppressed the
sequence index $k$, and since $\det(\gh)=1$, we omitted $dV_{\gh} =
d\xi$.  Finally, we remark that we chose to write this identity in
terms of $\uh_k$ rather than $v_k$ to avoid an even longer formula. It
is easy to translate the long range
estimates~\eqref{nov2e1}--\eqref{nov2e6} to $\uh_k$ from the fact
that $v_k(y)=\uh_k(\eps_k y)+\log \epsilon_k$.

We denote the integral on the left hand side of this identity by $I_0$, and the four integrals on the right-hand side by $I_1$, $I_2$, $I_3$ and $I_4$ respectively.  By~\eqref{oct31e3}, we obtain:
\[
    \frac12 \int_{\Omega} \xi^i\partial_i h_k e^{4\uh_k} = O(\epsilon_k).
\]
hence it follows from~\eqref{jan11e1} that:
\begin{equation}    \label{I0}
    I_0 = \alpha_k+O(L_k^{-3})
\end{equation}
By the expansions~\eqref{oct22e9} and~\eqref{oct22e3} of the derivatives of the metric $\gh$, and~\eqref{nov29e8} we get:
\begin{equation}    \label{I2}
    |I_2| \leq C \int_{B_{L_k}} \eps _k^3 |D^2v_k|\, |Dv_k|\, |y|^2 + \eps_k^2 |y|^2 |D^2v_k|^2
    + \eps_k^4 |Dv_k|\, |y| = O(\eps_k)
\end{equation}
and similarly, using~\eqref{nov29e1}, we see that:
\begin{equation}    \label{I3I4}
    |I_3| + |I_4| = O(\epsilon_k).
\end{equation}

It remains to compute $I_1$.  First, using~\eqref{oct31e3}, we can estimate the first term in $I_1$:
\[
    \frac12\int_{\partial\Omega} \xi^i \nu_i h_k e^{4\uh_k} = O(L_k^{-3}).
\]
Using this bound, and using the expansions~\eqref{oct22e0} and~\eqref{oct22e1} in the remaining terms, we can now reduce $I_1$ to:
\begin{align*}
    I_1 &= \int_{\partial\Omega} \left(
    -l_k \partial_\nu(\Delta\uh_k) \partial_\nu\uh_k + \Delta\uh_k \partial_\nu\uh_k
    + \nu_i\xi^m \Delta\uh_k \partial_{im}\uh_k - \frac12 l_k (\Delta\uh_k)^2 \right) \\
    & \qquad \qquad \qquad\qquad \qquad \qquad {} + O(\eps_k) \\
    & = \int_{\partial B_{L_k}} \left(
    -L_k \partial_\nu (\Delta v_k) \partial_\nu v_k + \partial_\nu (y\cdot\nabla v_k) \Delta v_k
    -\frac 12 L_k (\Delta v_k)^2 \right)
    + O(\eps_k).
\end{align*}
Using (\ref{nov2e1})-(\ref{nov2e6}) in the above, we get:
\begin{equation}    \label{I1}
    I_1 = \frac{\alpha_k^2}{16\pi^2} + O(L_k^{-1}).
\end{equation}
Combining~\eqref{I0}, \eqref{I2}, \eqref{I3I4} and~\eqref{I1}, we get:
\[
    \alpha_k + O(L_k^{-3}) = \frac{\alpha_k^2}{16\pi^2} + O(L_k^{-1}).
\]
which implies
\begin{equation}    \label{alphak}
    \alpha_k = 16 \pi^2 + O(L_k^{-1}).
\end{equation}
When substituting this into~\eqref{nov2e1}, we obtain:
\[
    v_k(y) + 2\log |y| = O(1), \quad |y|\ge L_k.
\]

The argument in the region $|y|\le L_k$ follows the one
in~\cite{linwei1} closely, hence we again only sketch the proof.
Without loss of generality, we assume that $h_k(q_k)\to 1$\footnote{Otherwise, we can
add a constant to $\hat u_k$.}.
 Let $U=U_{0,1,1}$ be the standard bubble in
$\R^4$:
\[
    U(y)= - \log \left(1 + \frac{|y|^2}{4\sqrt{3}}\right)
\]
It is easy to check that $U$ satisfies $\Delta^2U=2e^{4U}$ and it is well known that $v_k\to U$ in $C^4_{\text{loc}}(\mathbb R^4)$, see for example~\cite{druet1}.  Thus, for any fixed $A$ and all $k$ sufficiently large, we have:
\[
    |v_k(y)-U(y)|\le 1, \qquad \text{for $|y|\le A$.}
\]
Subtracting~\eqref{nov30e3} from its Euclidean counterpart, using~\eqref{mar19e1} to compare $d_{\gb}(y,z)$ and $|y-z|$ as well as their respective derivatives, and also using~\eqref{alphak} to compare the leading terms, we obtain:
\[
    |\Delta v_k(y) - \Delta U(y)| \le  C |y|^{-3}, \qquad \text{for $A<|y|< L_k$.}
\]
Now, letting $T(y)=C(1+|y|^{-1})$, and choosing $C$ large enough, we can guarantee that $\Delta T \le - |\Delta v_k - \Delta U|$
whence from the maximum principle $|v_k(y)-U(y)| \le C(1+|y|^{-1})$, on $A\le |y|\le L_k$.  Substituting $\xi=\eps_k y$, and using the definition of $v_k$, we obtain the version of~\eqref{MainEst} with $O(1)$ on the right-hand side, i.e., with $\tau=0$.

\section{A Sharper Estimate} \label{sec:03}

The main purpose of this section is to establish~\eqref{MainEst}. An important tool we use is the following lemma,
due to Lin-Wei~\cite{linwei1}:
\begin{lem} \label{nov5l1}
Let $U(y)=-\log (1+|y|^2/4\sqrt{3})$ be defined on $\R^4$.   Then $U$ satisfies
\[
    \Delta^2U=2e^{4U},\qquad U(0)=\max U = 0.
\]
Furthermore, any solution of the linearized problem:
\[
    \Delta^2\phi = 8e^{4U} \phi,  \qquad |\phi(y)|\le C(1+|y|)^{\tau}, \qquad \tau\in (0,1),
\]
is given by $\phi=\sum_{j=0}^4 c_j \psi_j$ where
\begin{eqnarray*}
    \psi_0 &=& \frac{1-|y|^2/4\sqrt{3}}{1+|y|^2/4\sqrt{3}} \\[1ex]
    \psi_j &=& \frac{y_j}{1+|y|^2/4\sqrt{3}}, \qquad j=1,..,4.
\end{eqnarray*}
\end{lem}

\begin{rem}
One immediate consequence of~Lemma~\ref{nov5l1} is that if in addition, $\phi$ satisfies
$\phi(0)=0$, and $\nabla\phi(0)=0$, then $\phi\equiv 0$.
 \end{rem}

Let $\rho_k=h_k(0)^{1/2}/4\sqrt3$ and consider the solution $U_k(y) = - \log (1+\rho_k |y|^2)$ of the equation
\begin{equation}    \label{Uk}
    \Delta^2U_k=2h_k(0)e^{4U_k}, \qquad U_k(0)=0,\qquad |\nabla U_k(0)|=0,
\end{equation}
on $\R^4$.  Letting $w_k=v_k-U_k$, then by the result of Section~\ref{sec:02}, we already know that $|w_k|\le C$ in
$B(0,\delta_1\eps_k^{-1})$.  Our goal in this section is to prove:
\begin{equation}    \label{nov21e1}
    |w_k(y)| \le C \eps_k^{\tau}|y|^{\tau},\qquad |y|\le \delta_1\eps_k^{-1}.
\end{equation}
for any $0<\tau<1$, which implies~\eqref{MainEst}.  Equivalently, if we let
\[
    \Lambda_k = \max_{\Omega_k} \frac{|w_k(y)|}{\eps_k^{\tau}(1+|y|)^{\tau}},
\]
then it suffices to show that $\Lambda_k$ is bounded on $\Omega_k=B(0,\delta_1\eps_k^{-1})$.
Suppose that $\Lambda_k\to \infty$, and let $y_k\in\Omega_k$ be the point where $\Lambda_k$ attains its maximum.  Now, either: (i) $y_k\to\infty$; or (ii) $|y_k|$ remains bounded at least along a subsequence, and hence a further subsequence, which without loss of generality we will assume is $y_k$ itself, converges to $y^*$.  We will show that in both cases a contradiction follows.

Define:
\[
    \bar w_k(y) =  \frac{w_k(y)}{\Lambda_k\eps_k^{\tau}(1+|y_k|)^{\tau}}.
\]
By the definition of $\Lambda_k$, we have
\begin{equation}    \label{nov12e2}
    |\bar w_k(y)| \le \left(\frac{1+|y|}{1+|y_k|}\right)^{\tau},
\end{equation}
and $\bar w_k(y_k) = \pm 1$.

Assume first that $y_k\to\infty$.  Since $|w_k(y)|\le C$ and $\Lambda_k\to\infty$, we clearly have $y_k=o(1)\eps_k^{-1}$.
>From the fundamental solution for $\Delta^2$, it is straightforward to get:
\begin{equation}
\begin{aligned}     \label{nov12e3}
    U_k(x) &= \frac 1{4\pi^2} \int_{\R^4} \log \frac{|y|}{|x-y|}\, h_k(0)\, e^{4U_k(y)}\, dy \\
    &= \frac 1{4\pi^2} \int_{\Omega_k} \log \frac{|z|}{|y-z|}\, h_k(0)\, e^{4U_k(z)}\,dz + O(\eps_k^4).
\end{aligned}
\end{equation}
Similarly, using the fundamental solution for $P_{\gb_k}$, we find:
\begin{equation}
\begin{aligned}  \label{jun29e1}
    v_k(y) &= \frac 1{4\pi^2}\int_{\Omega_k} \log\frac{|z|}{d_{\gb_k}(y,z)}\, h_k(\eps_kz)\, e^{4v_k(z)}\,dz + O(\eps_k|y|) \\
    &= \frac 1{4\pi^2} \int_{\Omega_k} \log\frac{|z|}{|y-z|}\, h_k(\eps_kz)\, e^{4v_k(z)}\, dz + O(\eps_k|y|),
\end{aligned}
\end{equation}
see~\eqref{oct31e5}.  Note that for the second equality, we used~\eqref{feb13e3} as well as
the decay rate of $v_k$.  Finally, we estimate the source term:
\begin{multline*}
    |h_k(\epsilon_kz)\,e^{4v_k(z)} - h_k(0)\,e^{4U_k}| \\
    = \left| h_k(\eps_k z)\,(e^{4v_k} - e^{4U_k}) + (h_k(\eps_kz)-h_k(0))\, e^{4U_k} \right| \\
     \le C (1+|z|)^{-8} |w_k(z)| + O(\eps_k)\,(1+|z|)^{-7} \\
     \le C\eps_k^{\tau} \Lambda_k(1+|z|)^{-8+\tau} + O(\eps_k)(1+|z|)^{-7}.
\end{multline*}
Substituting this in~\eqref{jun29e1} and combining with~\eqref{nov12e3} in the definition of $\bar w_k$,
we obtain:
\begin{multline*}
    \bar w_k(y_k) = \int_{\Omega_k}\log \frac{|z|}{|y-z|}
    \left( \frac{O(1)(1+|z|)^{-8+\tau}}{(1+|y_k|)^{\tau}} \right. \\
    \left. {} + \frac{O(\eps_k^{1-\tau})(1+|z|)^{-7}}{\Lambda_k(1+|y_k|)^{\tau}} \right)\, dz
    + o(1).
\end{multline*}
Since $y_k\to \infty$, it is now easy to see that the
right hand side is $o(1)$, which contradicts $\bar w_k=\pm 1$.

We now turn to the second case and assume without loss of generality that $y_k$ converges to $y^*$.  We will show that along a subsequence $\bar w_k$ converges.  This will be accomplished by estimating $P_{\gb_k}(U_k-v_k)$.  We start with:
\[
    P_{\gb} U_k = \Delta_{\gb}^2 U_k + \eps_k^2 \div_{\gb} \left( \left( \frac23
    \hat R (\eps_k y) \gb_{ij}(y) - 2 \hat R_{ij}(\eps_k y) \right)\, dU_k \right)
\]
where we have suppressed the subscript $k$ on the metric and curvature components, and where $\hat R$, $\hat R_{ij}$ are the scalar
and the Ricci curvatures of $\gh$.  In conformal normal coordinates, we have:
\begin{equation}    \label{1term}
    \Delta_{\gb}^2 U_k = \Delta^2 U_k \qquad \text{in $B(0,\delta_1 \epsilon_k^{-1})$.}
\end{equation}
Furthermore:
\begin{multline*}
    \partial_m \left( \gb^{mi}\left( \frac23 \hat R (\eps_k y) \gb_{ij} - 2\hat R_{ij}(\eps_ky) \right)\gb^{lj}\partial_l U_k \right) \\
    = \partial_m \gb^{mi} \left( \frac23 \hat R(\eps_k y) \gb_{ij} - 2\hat R_{ij}(\eps_k y) \right) \gb^{lj} \partial_l U_k \\
    {} + \gb^{mi} \partial_m \gb^{lj} \left( \frac23 \hat R (\eps_k y) \gb_{ij} - 2 \hat R_{ij}(\eps_k y) \right) \partial_l U_k \\
    {} + \eps_k \gb^{mi} \gb^{lj} \left( \frac 23 \partial_m \hat R(\eps_k y) \gb_{ij} + \frac23 \hat R(\eps_ky) \partial_m
    \gh_{ij}(\eps_ky) - 2\hat R_{ij,m}(\eps_ky) \right) \partial_l U_k\\
    {} + \gb^{mi} \gb^{lj} \left( \frac23 \hat R(\eps_ky) \gb_{ij} - 2\hat R_{ij}(\eps_k y) \right) \partial_{lm} U_k\\
    = A_1 + A_2 + A_3 + A_4.
\end{multline*}
Since $\partial_m\gb^{mi}(y)=\eps_k\partial_m\gh^{mi}(\eps_ky)=O(\eps_k^3|y|^2)$, $\hat
R(\eps_ky)=O(\eps_k^2|y|^2)$, and $\hat R_{ij}(\eps_ky)=O(\eps_k|y|)$, we easily get $A_1=O(\eps_k^2)$, and $A_2=O(\eps_k^2)$.  Furthermore, since in addition $\hat R_{ij,i}(0)=-\frac 12\hat R_{,j}(0)=0$, we also have $A_3=O(\eps_k^2)$, and $A_4=O(\eps_k^2)$.  Substituting this into the above equation and multiplying by $\eps_k^2$, we get:
\begin{equation}    \label{2term}
    \eps_k^2 \div_{\gb} \left( \left(\frac23 \hat R(\eps_k y) \gb_{ij}(y) - 2\hat R_{ij}(\eps_ky) \right)\, dU_k \right)
    = O(\eps_k^4).
\end{equation}
Combining~\eqref{1term} and~\eqref{2term}, we find:
\begin{equation} \label{nov5e1}
    P_{\gb} U_k = 2h_k(0)\,e^{4U_k} + O(\eps_k^4), \qquad |y|\le \delta_1\eps_k^{-1}.
\end{equation}
Combining~\eqref{nov5e1} with~\eqref{oct31e1}, we obtain:
\begin{equation}
    \label{nov5e2} P_{\gb} w_k = 8 h_k(\eps_k y)\,e^{4\xi^k} w_k + O(\eps_k) (1+|y|)^{-7} + O(\eps_k^4),
    \qquad |y| \le \delta_1 \eps_k^{-1},
\end{equation}
where $\xi^k$ is given by: $e^{4\xi^k}=\int_0^1e^{4tv_k+4(1-t)U_k}\,dt$.  Finally, this leads to the following equation for $\bar w_k$:
\begin{equation}
\label{nov5e3} P_{\breve g}\bar w_k=8h_k(\epsilon_ky)e^{4\xi^k}\bar
w_k+
\frac{O(\epsilon_k^{1-\tau})(1+|y|)^{-7}}{\Lambda_k(1+|y_k|)^{\tau}}+\frac{O(\epsilon_k^{4-\tau})}{\Lambda_k(1+|y_k|)^{\tau}}.
\end{equation}
Since $y_k\to y^*$, a subsequence of $\bar w_k$ converges to $w$ in $C^4(\R^4)$.  We will assume without loss of generality, as in Section~\ref{sec:02}, that $h_k(0)\to 1$.  It follows that the limit $w$ satisfies:
\[
    \begin{cases}
            \Delta^2w=8e^{4U}w,\\
            |w(y)| \le C(1+|y|)^{\tau},\\
            w(0)= |\nabla w(0)| =0.
    \end{cases}
\]
By the remark following Lemma~\ref{nov5l1}, we conclude that $w\equiv 0$, which contradicts $\bar w(y^*)=\pm 1$.  This concludes the proof of~\eqref{nov21e1}.

In Section~\ref{sec:04}, we will also need estimates on the derivatives $D^jw_k(y)$ for $|y|\le \eps_k^{-1}\delta_1/2$, $j=1,2,3$. By combining~\eqref{nov5e2} and~\eqref{nov21e1}, we have:
\[
    P_{\gb} w_k(y) = O(\eps_k^{\tau}) (1+|y|)^{-8+\tau} + O(\eps_k^4), \qquad  |y|< \frac{\delta_1}2\epsilon_k^{-1}.
\]
Fix $y$, and let $r=|y|$ and $f_k(z)=w_k(rz)$ for $1/2< |z| <2$.  Then $f_k(z)$ satisfies:
\begin{gather*}
    P_{\acute g} f_k(z)=O(\eps_k^{\tau}) (1+r)^{-4+\tau} + O(\eps_k^4 r^4), \qquad B_2\setminus B_{1/2}, \\
    f_k(z)=O(\eps_k^{\tau}r^{\tau}),\quad B_2\setminus B_{1/2},
 \end{gather*}
where $\acute g$ is the rescaled metric $r^{-2}\psi_*\gb$ and $\psi\colon z\mapsto rz$.  Standard elliptic theory for fourth order equations~\cite{browder} yields:
\[
    D^j f_k(z) = O(\eps_k^{\tau} r^{\tau}), \qquad |z|=1, \quad j=1,2,3.
\]
Hence, we conclude:
\begin{equation}\label{feb19e1}
    D^j w_k(y) = O(\eps_k^{\tau} |y|^{\tau-j}),  \qquad j=1,2,3,  \quad |y| \leq \eps_k^{-1}\delta_1/2.
\end{equation}

\section{The Vanishing Rate}
\label{sec:04}

The purpose of this section is to complete the proof of our Theorem by proving~\eqref{vrate}.  We first prove a Euclidean version in Subsection~\ref{sec:vrate}, and then translate this result to the original metric in Subsection~\ref{sec:original}.

\subsection{A Euclidean Version}        \label{sec:vrate}
The goal of this subsection is to prove~\eqref{nov22e8}.  This is accomplished in three steps:
\begin{enumerate}
 \item[(i)] In the first step, we derive an asymptotic expansion of the Pohozaev identity; see~\eqref{nov21e4}.
 \item[(ii)] In the second step, we express this identity in terms of the background Euclidean metric; see~\eqref{dec6e2}.
 \item[(iii)]  In the last step, we complete the proof of~\eqref{nov22e8}.
\end{enumerate}

\subsubsection{Step 1}
In this subsection we let $E_k=B(0,\eps_k^{-1/2})$, and we derive an asymptotic Pohozaev identity for $v_k$ on $E_k$. We multiply~\eqref{oct31e1} by $\partial_av_k$, $a=1,..,4$, integrate with respect to the Euclidean volume element $dy$, and estimate each of the resulting terms.  First by the $O(1)$ estimate and~\eqref{b}:
\begin{multline}  \label{nov21e3}
    \int_{E_k} 2 h_k(\eps_ky)\, e^{4v_k(y)}\, \partial_a v_k(y)\, dy \\
    = - \frac{\eps_k} 2\int_{E_k} e^{4v_k(y)}\, \partial_a h_k(\eps_ky)
    + \frac12 \int_{\partial E_k} h_k\, e^{4v_k}\, \nu_a \\
    = - \frac{\eps_k} 2\partial_a h_k(0) \int_{E_k} e^{4v_k(y)} + O(\eps_k^2).
\end{multline}
Next, integrating by parts, we have:
\begin{multline}    \label{Ek}
    \int_{E_k} \Delta_{\gb}^2v_k\,\partial_av_k \\
    = \int_{\partial E_k} \left(\gb^{ij}\partial_j(\Delta_{\gb}v_k)\partial_a v_k\,\nu_i
    -\gb^{ij}\Delta_{\gb}v_k\,\partial_{ia}v_k\,\nu_j
    + \frac 12(\Delta_{\gb}v_k)^2\,\nu_a \right)\\
    -\int_{E_k}(\Delta_{\gb}v_k \,\partial_{ia}\gb^{ij}\,\partial_j v_k + (\Delta_{\gb}v_k)\,
    \partial_a\gb^{ij}\,\partial_{ij}v_k),
\end{multline}
where we used:
\[
\partial_j(\breve g^{ij}\partial_{ia}v_k)=\partial_a(\Delta_{\breve g}v_k)
-\partial_{ia}\breve g^{ij}\partial_jv_k-\partial_a\breve
g^{ij}\partial_{ij}v_k
\]
We now estimate the two integrals over $E_k$ in~\eqref{Ek} above.  Using $\partial_{ia}\breve g^{ij}=O(\epsilon_k^3|y|)$, which is
implied by~\eqref{nov29e1}, and~\eqref{nov29e8}, we find:
\[
    \int_{E_k}\Delta_{\gb} v_k\, \partial_{ia}\gb^{ij}\, \partial_j\, v_k=O(\eps_k^2).
\]
Furthermore, for any $0<\sigma<1$, the second integral over $E_k$ in~\eqref{Ek} can be estimated as follows:
\begin{multline*}
    \int_{E_k}(\Delta_{\gb}v_k)\, \partial_a\gb^{ij}\,\partial_{ij}v_k \\
    =\eps_k^2\int_{E_k}(\Delta U_k + \Delta_{\gb}w_k)\, \left(-\frac
    23\Rh_{i(am)j}(0)\,y^m + O(\eps_k|y|^2) \right)\partial_{ij}v_k\\
    =-\frac 23\eps_k^2\int_{E_k}\Delta
    U_k\Rh_{i(am)j}(0)\,y^m\partial_{ij}U_k\,dy+O\bigl(\eps_k^{(3+\sigma)/2}\bigr)=O\bigl(\eps_k^{(3+\sigma)/2}\bigr),
    \end{multline*}
where we used~\eqref{feb19e1}, and the following estimate implied by~\eqref{oct22e9}:
\[
    \partial_a\gb^{ij}(y)=\eps_k\, \partial_a\gh^{ij}(\eps_ky)=-\frac
    23\Rh_{i(am)j}(0)\, \eps_k^2\, y^m + O(\eps_k^3|y|^2),
\]
as well as the antisymmetry of the curvature tensor, and $\Rh_{ij}(0)=0$.  Here, we use the customary round brackets notation to denote the symmetric part.  We will choose $\sigma\in(0,1)$ at the end of the argument of Section~\eqref{sec:vrate}.  Next, since $\Rh(\eps_k y)= O(\eps_k^2|y|^2)$, and $\Rh_{ij}(\eps_k y)=O(\eps_k |y|)$ and~\eqref{nov29e8}, we have:
\begin{multline*}
    \eps_k^2\int_{E_k}\partial_m\left( \gb^{mi}\left( \frac
    23\hat R(\eps_ky) \, \gb_{ij}-2\hat R_{ij}(\eps_ky) \right)\partial_l v_k\, \gb^{lj} \right)\partial_a v_k\\
    = \eps_k^2\int_{\partial E_k} \gb^{mi} \left(\frac
    23\Rh(\eps_ky)\,\gb_{ij}-2\Rh_{ij}(\eps_ky) \right)\, \partial_l v_k\, \gb^{lj}\, \partial_a v_k\, \nu_m\\
    - \eps_k^2 \int_{E_k} \gb^{mi}  \left(\frac
    23\Rh(\eps_ky) \, \gb_{ij}-2\Rh_{ij}(\eps_ky) \right)\partial_l v_k\,\gb^{lj}\, \partial_{am} v_k = O(\eps_k^2).
\end{multline*}
Finally, we estimate:
\[
    \eps_k^4 \int_{E_k} 2 \bh_k \partial_a v_k = O(\eps_k^{5/2}).
\]
Combining all the terms, we arrive at the following Pohozaev identity:
\begin{multline}    \label{nov21e4}
    \frac{\eps_k}2\, \partial_a h_k(0) \int_{E_k} e^{4v_k} \, dy+ O\bigl(\eps_k^{(3+\sigma)/2}\bigr) \\
    = \int_{\partial E_k} \left( -\gb^{ij}\,\partial_j(\Delta_{\gb}v_k)\,\partial_a v_k\, \nu_i
    + \gb^{ij}\,\Delta_{\gb}v_k\,\partial_{ia}v_k\,\nu_j - \frac12 (\Delta_{\gb}v_k)^2\, \nu_a \right).
\end{multline}

\subsubsection{Step 2}
In this second step, we rewrite~\eqref{nov21e4} in terms of the Euclidean $\Delta v_k$ rather than $\Delta_{\gb} v_k$.  We begin by substituting:
\[
    (\gb^{ij}(y)-\delta_{ij})\,\nu_i = -\frac13 \eps_k^2 \Rh_{ilmj}(0)\,y_m\,y_l\,\frac{y_i}{|y|} + O(\eps_k^3|y|^3) = O(\eps_k^3|y|^3).
\]
into~\eqref{nov21e4} to get:
\begin{multline} \label{nov21e5}
    \frac{\eps_k}2\,\partial_a h_k(0) \int_{E_k} e^{4v_k}\,dy + O\bigl(\eps_k^{(3+\sigma)/2}\bigr) \\
    = \int_{\partial E_k} \left( - \partial_i(\Delta_{\gb} v_k)\,\partial_a v_k \, \nu_i + \Delta_{\gb}v_k\, \partial_{ia} v_k\, \nu_i
    -\frac 12 (\Delta_{\gb}v_k)^2\,\nu_a \right).
\end{multline}
A straightforward computation leads to:
\begin{multline}    \label{jun30e1}
    \partial_i(\Delta_{\gb}v_k)-\partial_i(\Delta v_k)\\
    = \partial_{im}\gb^{ml}\,\partial_l v_k
    + \partial_m \gb^{ml}\,\partial_{il}v_k + \partial_i\gb^{ml}\, \partial_{ml}v_k + (\gb^{ml}-\delta_{ml})\, \partial_{iml}v_k.
\end{multline}
In view of~\eqref{nov29e1} and~\eqref{oct22e3}, we have, for $|y|=\eps_k^{-1/2}$:
\begin{gather}
    \label{jun30e2a}
    \partial_{im}\gb^{ml}\,\partial_l v_k = O(\eps_k^3),\qquad
    \partial_m \gb^{ml}\, \partial_{il}v_k = O(\eps_k^3),\\
    \label{jun30e2}
    \partial_i\gb^{ml}\, \partial_{ml}v_k = \partial_i\gb^{ml}\, (\partial_{ml}U_k
    + O(\eps_k^{\sigma}r^{\sigma-2})) = O\bigl(\eps_k^{(5+\sigma)/2}\bigr),
\end{gather}
where to derive~\eqref{jun30e2}, we also used the following consequence of~\eqref{oct22e9}:
\[
    \partial_i \gb^{ml} = -\frac23 \Rh_{m(ia)l}(0) \eps_k^2 y_a + O(\eps_k^3|y|^2),
\]
as well as the anti-symmetry of the curvature tensor and $\Rh_{ij}(0)=0$.
Next, on $|y|=\eps_k^{-1/2}$, we have
\begin{multline} \label{jun30e3}
    (\gb^{ml}-\delta^{ml})\, \partial_{iml}v_k \\
    = \left(-\frac13\eps_k^2 \Rh_{mabl}(0)\,y_a\,y_b + O(\eps_k^3r^3) \right)
        (\partial_{iml}U_k + O(\eps_k^{\sigma}r^{\sigma-3})) \\
    = -\frac13\eps_k^2\Rh_{mabl}(0)\,y_a\,y_b\,\partial_{iml}U_k + O\bigl(\eps_k^{(5+\sigma)/2}\bigr)
    = O\bigl(\eps_k^{(5+\sigma)/2}\bigr),
\end{multline}
where we have used the following expansion, valid for any radial function $f(r)$:
\begin{multline*}
    \partial_{iml} f(r) = \left( f'''(r)-\frac{f''(r)}{r} + \frac{f'(r)}{r^2} \right) \frac{y_l\,y_m\,y_i}{r^3} \\
    {}  + \left( f''(r)-\frac{f'(r)}r \right) \frac{(\delta_{il}\,y_m + y_l\, \delta_{im}) r^2 - 2y_l\,y_m\,y_i}{r^4} \\
    {} + (f''(r)-f'(r))\frac{\delta_{ml}\,y_i}{r^2},
\end{multline*}
as well as the anti-symmetry of $\Rh_{abcd}$ and $\Rh_{ij}(0)=0$.
It now follows from~\eqref{jun30e1}, \eqref{jun30e2} and~\eqref{jun30e3} that the following holds:
\begin{equation} \label{jun30e4}
    -\int_{\partial E_k} \partial_i(\Delta_{\gb}v_k)\,\partial_a v_k\,\nu_i
    = -\int_{\partial E_k} \partial_i(\Delta v_k)\, \partial_a v_k\, \nu_i + O\bigl(\eps_k^{(3+\sigma)/2}\bigr).
\end{equation}
Next, on $|y|=\epsilon_k^{-1/2}$, we have:
\begin{multline*}
    \Delta_{\gb}v_k=\partial_i\gb^{ij}\,\partial_jv_k + \gb^{ij}\,\partial_{ij}v_k
    = O(\eps_k^3)r + \Delta v_k+ (\gb^{ij}-\delta_{ij})\,\partial_{ij}v_k\\
    =O(\eps_k^{\frac 52})+\Delta v_k- \left( \frac{\eps_k^2}3 \Rh_{iabj}(0)\,y_a\,y_b + O(\eps_k^3r^3) \right)
    \bigl( \partial_{ij}U_k + O(\eps_k^{\sigma}r^{\sigma-2}) \bigr)\\
    = \Delta v_k+O\bigl(\eps_k^{2+\sigma/2}\bigr),
\end{multline*}
from which it follows:
\begin{equation} \label{jun30e5}
    \int_{\partial E_k} \Delta_{\gb}v_k\,\partial_{ia} v_k\,\nu_i
    = \int_{\partial E_k}\Delta v_k\,\partial_{ia} v_k\, \nu_i\, dS + O\bigl(\eps_k^{(3+\sigma)/2}\bigr).
\end{equation}
Similarly:
\begin{equation} \label{jun30e6}
    -\frac12 \int_{\partial E_k}(\Delta_{\gb}v_k)^2\, \nu_a
    = -\frac12 \int_{\partial E_k}(\Delta v_k)^2\,\nu_a + O\bigl(\eps_k^{(3+\sigma)/2}\bigr).
\end{equation}
Substituting~\eqref{jun30e4}, \eqref{jun30e5} and~\eqref{jun30e6} into the Pohozaev Identity~\eqref{nov21e5}, we conclude:
\begin{multline} \label{dec6e2}
    \frac{\eps_k}2\partial_a h_k(0) \int_{E_k} e^{4v_k}\, dy + O(\eps_k^{(3+\sigma)/2}) \\
    = \int_{\partial E_k} \left(-\partial_i(\Delta v_k)\,\partial_a v_k\,\nu_i + \Delta v_k\, \partial_{ia} v_k\, \nu_i
    - \frac12(\Delta v_k)^2\,\nu_a \right)\,dS.
\end{multline}

\subsubsection{Step 3} In this subsection, we first aim to replace $v_k$ by $U_k$ in the Pohozaev identity~\eqref{dec6e2}, after which many of the terms will simplify thanks to the radial symmetry of $U_k$, leading to the Euclidean version of the vanishing rate~\eqref{nov22e8}.
Recall the definition of $\phih_k$~\eqref{feb15e1}, from which we have:
\[
    D^j\phih_k(\eps_ky)=O(\eps_k^j),\qquad j=1,2,3.
\]
For $|y|=\eps_k^{-1/2}$, we cut $B(0,\delta_1\eps_k^{-1})$ into three subdomains $B(0,\delta_1\eps_k^{-1})=\Omega_1\cup\Omega_2\cup\Omega_3$ as in Section~\ref{sec:02}, page~\pageref{omega}, and use the representation~\eqref{oct31e5}.
Using standard estimates over $\Omega_{2}\cup \Omega_{3}$ and~\eqref{jun27e10} over $\Omega_1$, we find:
\begin{multline*}
    \partial_a v_k(y) \\
    = -\frac 1{4\pi^2} \int_{\Omega_{1}}\partial_a \bigl(\log
    d_{\gb}(z,y) \bigr) \, \hb_k(z)\, e^{4v_k(z)} + \eps_k\, \partial_a\phih_k(0)
    {} + O(\eps_k^2|y|)+O(|y|^{-5}) \\
    = -\frac1{4\pi^2}\int_{\Omega_{1}} \frac{y_a-z_a}{|y-z|^2}\, h_k(\eps_kz)\, e^{4v_k(z)}
    + \eps_k\, \partial_a\phih_k(0) + O(\eps_k^{3/2}),
\end{multline*}
where we have omitted the standard volume element $dz$.  Similarly,
\begin{gather*}
    \Delta v_k(y)=-\frac 1{2\pi^2}\int_{\Omega_{1}}
    \frac 1{|y-z|^2}\, h_k(\eps_kz)\, e^{4v_k(z)} + O(\eps_k^2), \\[1ex]
    \partial_{ij}v_k(y)=-\frac 1{4\pi^2}\int_{\Omega_{1}}
    \frac{\delta_{ij}|y-z|^2-2(y_i-z_i)(y_j-z_j)}{|y-z|^4}\, h_k(\eps_kz)\, e^{4v_k}
    +O(\eps_k^2). \\[1ex]
    \partial_i(\Delta v_k(y))=\frac 1{\pi^2}\int_{\Omega_{1}}
    \frac{y_i-z_i}{|y-z|^4}\, h_k(\eps_kz)\, e^{4v_k} + O(\eps_k^{5/2}).
\end{gather*}
We also have
\begin{align*}
    h_k(\eps_kz)\, e^{4v_k(z)}&= \bigl(h_k(0)+ \eps_k\, \partial_j h_k(0)\, z_j+O(\eps_k^2|z|^2)
     \bigr)e^{4U_k+O(\eps_k^{\sigma}|z|^{\sigma})}\\
    &=h_k(0)\, e^{4U_k(z)}+O(\eps_k^{\sigma})(1+|z|)^{-8+\sigma}.
\end{align*}
Substituting this in the derivatives of $v_k$ we have:
\begin{gather*}
    \partial_av_k(y)=-\frac1{4\pi^2}\int_{\Omega_{1}}\frac{y_a-z_a}{|y-z|^2}\,h_k(0)\,e^{4U_k(z)}
    +\eps_k\,\partial_a\phih_k(0)+O(\eps_k^{\sigma+1/2})\\[1ex]
    \Delta v_k(y)=-\frac 1{2\pi^2}\int_{\Omega_{1}} \frac1{|y-z|^2}\,h_k(0)\,e^{4U_k(z)}+O(\eps_k^{\sigma+1}),\\[1ex]
    \partial_{ij}v_k(y)=-\frac 1{4\pi^2}\int_{\Omega_{1}}\frac{\delta_{ij}\,|y-z|^2-2(y_i-z_i)(y_j-z_j)}{|y-z|^4}\,h_k(0)\,e^{4U_k}
    +O(\eps_k^{\sigma+1})\\[1ex]
    \partial_i(\Delta v_k(y))=\frac 1{\pi^2}\int_{\Omega_{1}}
    \frac{y_i-z_i}{|y-z|^4}\,h_k(0)\,e^{4U_k}+O(\eps_k^{\sigma+3/2}).
\end{gather*}
We perform a similar computation for $U_k$ and take the difference, leading to the following estimates for
$|y|=\eps_k^{-1/2}$:
\begin{align*}
    \partial_i(\Delta v_k(y))&=\partial_i(\Delta
    U_k)(y)+O(\eps_k^{\sigma+3/2}),\\
    \Delta v_k(y)&=\Delta U_k(y)+O(\eps_k^{\sigma+1}),\\
    \partial_{ij}v_k(y)&=\partial_{ij}U_k(y)+O(\eps_k^{\sigma+1})\\
    \partial_av_k(y)&=\partial_aU_k(y)+\eps_k\,\partial_a\phih_k(0)+O(\eps_k^{\sigma+1/2}).
\end{align*}
Substituting these estimates into the Pohozaev Identity~\eqref{dec6e2}, we obtain:
\begin{multline}    \label{PUk}
    \frac{\eps_k}2\partial_ah_k(0)\int_{E_k}e^{4U_k} + O(\eps_k^{1+\sigma/2})\\
    =\int_{\partial E_k}\left(-\partial_{\nu}(\Delta
    U_k)(\partial_aU_k+\eps_k\,\partial_a\phih_k(0)) + \Delta
    U_k\,\partial_{ia}U_k\,\nu_i-\frac 12(\Delta U_k)^2\,\nu_a\right).
\end{multline}
The symmetry of $U_k$ implies:
\[
    \int_{\partial E_k}\left(-\partial_{\nu}(\Delta U_k)\,\partial_aU_k
    +\Delta U_k\,\partial_{ia}U_k\,\nu_i-\frac 12(\Delta U_k)^2\,\nu_a\right)=0.
\]
In view of the equation~\eqref{Uk}, we also have
\[
    \int_{\partial E_k}\partial_{\nu}(\Delta
    U_k)\,dS = 2h_k(0)\int_{E_k}e^{4U_k}.
\]
Substituting into~\eqref{PUk}, we obtain
\[
    \partial_ah_k(0)+4h_k(0)\partial_a
    \phih_k(0)=O(\eps_k^{\sigma/2})+O(\eps_k^{\sigma-1/2}),\quad a=1,2,3,4.
\]
Now, if we choose $\tau/2 + 1/2<\sigma<1$, then $O(\eps_k^{\sigma/2})+O(\eps_k^{\sigma-1/2})=O(\eps_k^{\tau/2})$, so that we can conclude:
\begin{equation}    \label{nov22e8}
    \left| \frac{\nabla h_k(0)}{h_k(0)}+4\nabla\phih_k(0) \right|= O(\eps_k^{\tau/2}).
\end{equation}

\subsection{The Vanishing Rate in $g$}  \label{sec:original}
In this subsection, we verify that~\eqref{nov22e8} leads to~\eqref{vrate}.  For simplicity, we assume without loss of generality that the cut-off function $\chi$ is supported in $B(q_{ik},2\delta)$ where $\delta$ is small enough to guarantee that $B(q_{ik},2\delta)$ are mutually
disjoint. Indeed, this can be done since the left hand side of~\eqref{vrate} is invariant under any change of cut-off function $\chi$.  Under this choice of cut-off, all the terms $\nabla_1 G(q_{ik},q_{jk})$, $j\ne i$, reduce to $\nabla_1 \beta(q_{ik},q_{jk})$, so that it now suffices to show that:
\begin{equation}    \label{mar25e1}
    64\pi^2\sum_{j=1}^N\nabla_1\beta(q_{ik},q_{jk}) -4\nabla\phi_k(q_{ik}) =-\frac{\nabla h_k(q_{ik})}{h_k(q_{ik})}+O(\eps_k^{\tau/2}).
\end{equation}
Indeed, $e^{-\tau u_k(q_{ik})/2}=O(\eps_k^{\tau/2})$ since $|u_k(q_{ik})-u_k(q_{jk})|\le C$ for $i\neq j$, and furthermore
$\nabla_1\beta(q_{ik},q_{ik})=\nabla_2\beta(q_{ik},q_{ik})$ since $\beta(x,y)=\beta(y,x)$ by Lemma \ref{jan18l1} in the Appendix.  The remainder of this section is devoted to verifying~\eqref{mar25e1}.

Taking the derivative with respect to $\xi$ in (\ref{uhat}) and
evaluating at $q_{ik}$, we have:
\begin{equation}    \label{uk}
    \nabla \hat
    u_k(q_{ik})=2\int_M\nabla_1H(q_{ik},\eta)h_k(\eta)e^{4\uh_k(\eta)}dV_{\gh}(\eta)+\nabla \phih_k(q_{ik}).
\end{equation}
where $H(\xi,\eta)=-(1/8\pi^2)\, \chi(r)\, \log d_{\gh_k}(\xi,\eta)$.  Similarly, for $u_k$ we have:
\begin{multline}    \label{ukhat}
    \nabla u_k(q_{ik})=-\frac{1}{8\pi^2}\int_M\nabla_{1} \bigl( \chi(r)\log d_g(q_{ik},\eta) \bigr) 2h_k(\eta)
    e^{ 4\uh_k(\eta)}dV_{\hat g}(\eta)\\
    {} + \int_M\nabla_1\beta(q_{ik},\eta)2h_k(\eta)e^{4\hat
    u_k(\eta)}dV_{\hat g}(\eta)-\nabla \phi_k(q_{ik}),
\end{multline}
where we used $e^{4\uh_k}\,dV_{\gh_k} = e^{4u_k}\, dV_g$.
Let $H_0(\xi,\eta) = (1/8\pi^2) \chi(r)\log d_g(\xi,\eta)$, then we claim that:
\begin{equation}    \label{HH}
    |\nabla_1H(q_{ik},\eta)-\nabla_1H_0(q_{ik},\eta)|\le C d_g(\xi,\eta).
\end{equation}
Indeed, recall that $\wh(q_{ik})=0$ and $\nabla \wh_k(q_{ik})=0$.  Thus, if fix $\xi=q_{ik}$, and we let:
\[
    f(\eta) = \log d_g(\xi,\eta) - \log d_{\gh}(\xi,\eta),
\]
then $\wh_k(\eta) = O(|\xi-\eta|^2)$, and therefore $|f(\eta)|\leq C(|\xi-\eta|^2)$.  It follows that $\nabla f(\xi) = 0$.  Now, by Appendix~\ref{appGreen},
$P_gf(\eta)$ is a bounded function, hence by elliptic theory, $\nabla^2 f(\eta)$ is bounded.  We conclude that $|\nabla f(\eta)|\leq C|\xi-\eta|$ from which~\eqref{HH} follows.  Next, since we have:
\[
    \nabla u_k(q_{ik})=\nabla \hat u_k(q_{ik})+\nabla \hat
    w_k(q_{ik})=\nabla \hat u_k(q_{ik}),
\]
it follows, by taking the difference of~\eqref{uk} and~\eqref{ukhat}, that:
\begin{equation}    \label{Dphik}
    0 = -\nabla \phih_k(q_{ik})+\int_M\nabla_1\beta(q_{ik},\eta)2h_k(\eta)e^{4\hat
    u_k(\eta)}dV_{\hat g}(\eta)-\nabla \phi_k(q_{ik})+O(\eps_k)
\end{equation}
Furthermore, we claim that:
\begin{equation}    \label{betak}
    \int_M\nabla_1\beta(q_{ik},\eta)2h_k(\eta)e^{4\hat
    u_k(\eta)}dV_{\hat g}(\eta) = 16\pi^2 \sum_{j=1}^N \nabla_1\beta(q_{ik},q_{jk}) + O(\eps_k^\tau).
\end{equation}
Indeed, observe that:
\begin{multline*}
    \int_M\nabla_1\beta(q_{ik},\eta)2h_k(\eta)e^{4\hat
    u_k(\eta)}dV_{\hat g} \\
    =\sum_{j=1}^N\int_{B(q_{jk},\delta)}\nabla_1\beta(q_{ik},\eta) 2 h_k(\eta) e^{4\uh_k(\eta)}\, dV_{\hat
    g}(\eta)+O(\eps_k^4)
\end{multline*}
since, by~\eqref{MainEst}, $e^{4\hat u_k}=O(\eps_k^4)$ on $M\setminus \cup_{j=1}^N B(q_{jk},\delta)$.
In addition, for each $j=1,..,N$:
\begin{multline*}
    \int_{B(q_{jk},\delta)}\nabla_1\beta(q_{ik},\eta)2h_k(\eta)e^{4\hat
    u_k(\eta)}dV_{\hat g}(\eta)\\
    =\nabla_1\beta(q_{ik},q_{jk})\int_{B(q_{jk},\delta)}2h_ke^{4\hat u_k}dV_{\hat
    g}+\int_{B(q_{jk},\delta)}O(|\eta-q_{jk}|)2h_ke^{4\hat u_k}dV_{\hat g}\\
    =16\pi^2+O(\epsilon_k^{\tau}),
\end{multline*}
where we used~\eqref{MainEst} to estimate the first integral and a standard rescaling to estimate the second one.
By combining~\eqref{nov22e8} with~\eqref{Dphik} and~\eqref{betak}, it follows that \eqref{mar25e1} holds.  This completes the proof of the Theorem.

\begin{rem}
Integrating~\eqref{jan22e1}, and using~\eqref{MainEst} on the right hand side, one easily obtains:
\[
    \int_M b_k\, dV_g = 8\pi^2 N + O(\eps_k^{\tau}).
\]
\end{rem}

\appendix

\section{The Green's function for $P_g$}    \label{appGreen}

Let $G$ denote the Green's function of $P_g$ as in~\eqref{green}:
\begin{equation}    \label{green1}
    f(\xi)-\bar f_g=\int_MG(\xi,\eta) P_gf(\eta) \, dV_g(\eta),
\end{equation}
where $\bar f_g=\Vol_g(M)^{-1}\int_M f\, dV_g$ is the mean value
of $f$ with respect to $g$. Clearly, $G$ is determined up to an arbitrary
function of $\xi$ which we can fix by imposing the condition:
\begin{equation}    \label{intG=0}
        \int_M G(\xi,\eta)\, dV_g(\eta)=0
\end{equation}

The purpose of this appendix is to prove the following lemma, which is an improvement on
a result of Chang-Yang~\cite{changyang1}:

\begin{lem} \label{jan18l1}
Let $(M,g)$ be a compact closed $4$-dimensional manifold, and suppose $\Ker(P_g)=\{\text{constants}\}$.
Then the Green's function $G(\xi,\eta)$
with respect to $P_g$ can be written as:
\begin{equation}    \label{eqG}
    G(\xi,\eta) = - \frac{1}{8\pi^2} \chi(r) \log r  + \beta(\xi,\eta)
\end{equation}
where  $r=d_g(\xi,\eta)$ is the geodesic distance between $\xi$
and $\eta$, $\chi(r)$ is a cut off function that is $1$ on a
neighborhood of $\xi$ and vanishes outside
$B\bigl(\xi,\delta(\xi)/10\bigr)$, and $\delta(\xi)$ is the
injectivity radius of $(M,g)$ at $\xi$.  Furthermore,
$\beta(\xi,\eta)\in W^{4,q}(M\times M)$, for any $1<q<\infty$ and
satisfies:
\[
    \|\beta(\xi,\cdot)\|_{W^{4,q}(M)} \le C, \qquad \text{uniformly in $\xi\in M$}
\]
for some constant $C=C(g,q)$.  The principal part of $G$ satisfies weakly:
\begin{equation}\label{jun27e1}
    P_{g,\eta} \left( - \frac 1{8\pi^2} \chi(r) \log d_g(\xi,\eta) \right) = \delta_{\xi} +  E(\xi,\eta)
\end{equation}
where $E$ is a bounded function.  Finally, we have $G(\xi,\eta)=G(\eta,\xi)$.
\end{lem}

\begin{proof}[Proof of Lemma~\ref{jan18l1}:]
The weak form of~\eqref{green1} is:
\begin{equation}    \label{weakgreen}
    P_{g,\eta}G(\xi,\eta)=\delta_{\xi}-\frac 1{\Vol_g(M)}.
\end{equation}
Let $\gbar=e^{2w}g$ be a metric conformal to $g$ such that in a neighborhood of $\xi$ the normal coordinates with respect to $\gbar$ are conformal normal coordinates, i.e., $\det(\gbar)=1$.  Similar to~\eqref{weakgreen}, we have for the Green's function $G_1$ of $P_{\gbar}$:
\[
    P_{\gbar,\eta}G_1(\xi,\eta) = \delta_{\xi} - \frac 1{\Vol_{\gbar}(M)}.
\]
By the conformal covariance of $P_{\gbar}$ this is equivalent to:
\[
    e^{-4w} P_{g,\eta} G_1 (\xi,\eta) = \delta_{\xi} - \frac{1}{\Vol_{\gbar}(M)}.
\]
Therefore since $w(\xi)=0$, we have:
\[
    P_{g,\eta}(G_1(\xi,\eta) - G(\xi,\eta)) = \frac1{\Vol_g(M)} - \frac{e^{4w(\eta)}}{\Vol_{\gbar}(M)}.
\]
Note that the integral of the right hand side is zero, hence there exists $F\in C^5(M\times M\setminus D)$ such that:
\[
    P_{g,\eta}F(\xi,\eta) = \frac 1{\Vol_g(M)} - \frac{e^{4w(\eta)}}{\Vol_{\gbar}(M)},
    \qquad \|\nabla_{\xi}F_1(\xi,\cdot)\|_{C^4(\cdot)}\le C.
\]
In particular, $G_1 = G + F$, and the principal part of $G$ and $G_1$ differ by a bounded function.  In the following, we focus on $G_1$.

In conformal normal coordinates the point $\xi$ corresponds to $0$.  We will identify $y\in T_\xi M$ with its image under the exponential map.  We first write $G_1$ as in~\eqref{eqG}:
\begin{equation}    \label{G1}
    G_1(0,y) = - \frac 1{8\pi^2} \chi(\rb) \log \rb + \beta_1(0,y),
\end{equation}
where $\rb=|y|=d_{\gbar}(0,y)$, and $\chi$ is a cut-off function.  We will show that the principal part $H = - (1/{8\pi^2}) \chi(\rb) \log \rb$ of $G_1$ satisfies weakly:
\begin{equation}    \label{eqH}
    P_{\gbar} H(0,y) = \delta_0 + E_1(0,y).
\end{equation}
where $E_1$ is a bounded function.  Since $P_{\gbar} \beta_1(0,y) = -E_1(0,y)$, this will imply by elliptic theory that $\beta_1(0,\cdot) \in W^{4,q}(M)\subset C^{3,\alpha}(M)$ for any $1<q<\infty$. Here we assume that $\chi\equiv 1$ in $B(0,\delta)$, $\chi\equiv 0$ in $M\setminus B(\xi,2\delta)$, and $\det(\gbar)\equiv 1$ in $B(0,2\delta)$.

Observe that since $H(0,y)$ is radial and $\chi$ is supported in a small neighborhood of $0$ where $\det(\gbar)=1$, we have $\Delta_{\gbar}^2 H(0,y) = \Delta^2 H(0,y)$. Thus, it follows that for any smooth function $\phi$:
\[
    \phi(0) = - \int_M\Delta_{\gbar}^2 H(0,y) \phi(y)\, dV_{\gbar}(y) + \int_M H(0,y) \Delta_{\gbar}^2 \phi(y)\, dV_{\gbar}(y).
\]
Clearly:
\[
    \int_M \Delta_{\gbar}^2 H(0,y) \phi(y)\, dV_{\gbar}(y)
    = \int_{B_{2\delta}\setminus B_{\delta}} \Delta_{\gbar}^2 H(0,y) \phi(y)\, dV_{\gbar}(y)
\]
and $\Delta_{\gbar}^2 H(0,y)$ is a bounded smooth function on $B_{2\delta}\setminus B_{\delta}$.  Hence weakly, we have:
\begin{equation}    \label{estH}
    \Delta_{\gbar}^2 H(0,y) = \delta_0 + \text{a bounded function}.
\end{equation}
Also, $(P_{\gbar}-\Delta_{\gbar}^2)H(0,y) = A \chi+ \text{a
bounded function}$, where:
\begin{eqnarray*}
    A &=& - \partial_m \left( \gbar^{mi} \left( \frac23 \Rb(y) \gbar_{ij}-2 \Rb_{ij}(y) \right) \gbar^{lj} \partial_l
    \left( \frac1{8\pi^2} \log \rb \right) \right) \\
    &=& - \frac 1{8\pi^2} \partial_m \left( \gbar^{mi} \left( \frac23 \Rb(y) \gbar_{ij} - 2\Rb_{ij}(y) \right)
    \gbar^{lj} \frac{y_l}{\rb^2} \right) \\
    &=& - \frac 1{8\pi^2}\partial_m \left( \left( \frac23 \Rb(y) \delta_{mj} - 2\Rb_{ij}(y) \gbar^{mi}\right)
    \gbar^{lj} \frac{y_l}{\rb^2}\right) \\
    &=& - \frac 1{8\pi^2} \left( \left( \frac23 \partial_m \Rb(y)\delta_{mj} -
    2\partial_m \bigl( \Rb_{ij}(y) \gbar^{mi}\bigr)\right)\gbar^{lj} \frac{y_l}{\rb^2}
    \right. \\ && \qquad\qquad + \left.
    \left( \frac23 \Rb(y) \delta_{mj} -
    2\Rb_{ij}(y) \gbar^{mi} \right) \partial_m \gbar^{lj} \frac{y_l}{\rb^2}
    \right. \\ && \qquad\qquad + \left.
    \left(\frac 23\bar R(y)\delta_{mj}-2\bar R_{ij}(y)\gbar^{mi} \right)\bar
    g^{lj}\frac{\delta_{lm}\rb^2-2y_ly_m}{\rb^4} \right)\\
    &=& - \frac 1{8\pi^2}(A_1 + A_2 + A_3).
\end{eqnarray*}
Here $\Rb$, $\Rb_{ij}$, and $\Rb_{ijkl}$ are used to denote the
scalar, Ricci, and Riemann curvatures of the metric $\gbar$. The
properties of conformal normal coordinates we used are listed in
Appendix~\ref{appcoords}.  We estimate each of the three terms in
the previous equation separately.
\begin{eqnarray*}
    A_1 &=& \left( \frac23 \partial_m \Rb(y) \delta_{mj} - 2\Rb_{ij,m}(y) \gbar^{mi}
    - 2\Rb_{ij}(y) \partial_m \gbar^{mi}\right) \gbar^{lj} \frac{y_l}{\rb^2} \\
    &=& \left( \frac23 \partial_m \Rb(y) \delta_{mj} - 2\Rb_{ij,i}(y) - 2\Rb_{ij}(y) \partial_m \gbar^{mi}\right)
    \frac{y_j}{\rb^2} + O(\rb) \\
    &=& - \frac13 \Rb_{,ij}(0) y_i y_j \rb^{-2} + O(\rb).
\end{eqnarray*}
In the second equality above, we used $\gbar^{mi} = \delta_{mi} + O(\rb^2)$ and $\gbar^{lj}=\delta_{lj} + O(\rb^2)$.
In the third equality, we used $\partial_m\Rb(0)=0$ and $2\Rb_{ij,im}(0) = \Rb_{,jm}(0)$, as well as $\partial_m \gbar^{mi}= O(\rb^2)$.
It is easy to see that $A_2=O(\rb)$ since $\Rb_{ij}(y)=O(\rb)$, $\Rb(y)=O(\rb^2)$, $\partial_m \gbar^{mi} = O(\rb)$.  Finally, we have
\begin{eqnarray*}
    A_3 &=& \left( \frac23 \Rb(y) \delta_{mj} - 2\Rb_{ij}(y) \gbar^{mi} \right) \gbar^{lj}\, \frac{\delta_{lm}\rb^2-2 y_l y_m}{\rb^4} \\
    &=& \left( \frac23 \Rb(y) \gbar^{lm} - 2\Rb_{ij}(y) \gbar^{mi} \gbar^{lj} \right) \, \frac{\delta_{lm}\rb^2-2 y_l y_m}{\rb^4} \\
&=& \left( \frac43 \Rb_{,ij}(0) y_i y_j \rb^2 - \frac43 \Rb(y) \rb^2 - 2 \Rb_{ij}(y) \delta_{mi} \delta_{mj} \rb^2
    \right. \\
    && \qquad\qquad {} + 4\Rb_{ij}(y) \delta_{mi} \delta_{lj} y_l y_m \biggr)  \rb^{-4} + O(\rb)\\
    &=& \left( -\frac13 \Rb_{,ij}(0) y_i y_j \rb^2 + 2\Rb_{ml,ab}(0) y_a y_b y_l y_m \right)  \rb^{-4} + O(\rb).
\end{eqnarray*}
In the third equality we used $\gbar^{ij}=\delta_{ij}+O(\rb^2)$, $\Rb(y)=O(\rb^2)$, $\Rb_{ij}(y)=O(\rb)$, and in the fourth
equality we used $\Rb_{ij,l}(0) y^i y^j y^l=0$. Combining the estimates for $A_1$, $A_2$ and $A_3$, we obtain:
\[
    A = \frac1{8\pi^2} \left( \frac23 \Rb_{,ij}(0) y_i y_j \rb^{-2} - 2\Rb_{ij,lm}(0) y_i y_j y_l y_m \rb^{-4} \right) + O(\rb).
\]
Finally, this last estimate, together with~\eqref{estH} yields~\eqref{eqH} as claimed earlier.

Rewriting~\eqref{eqH} in integral form, we get:
\begin{equation}\label{feb12e2}
    \phi(\xi) = \int_M H(0,y) P_{\gbar}\phi(y)\, dV_{\gbar}(y) - \int_M E_1(0,y) \phi(y)\, dV_{\gbar}(y).
\end{equation}
Note that $E_1(0,\cdot)$ is supported in $B(\xi,2\delta)\setminus
B(\xi,\delta)$.  In this small neighborhood, let $\eta =
\exp_\xi(y)$, where $\exp_\xi\colon T_\xi(M)\to M$ is the
exponential map with respect to the metric $g$. Thus, recalling
that the principal parts of $G$ and $G_1$ differ only by a bounded
function, we may rewrite~\eqref{feb12e2} as:
\begin{equation}    \label{feb12e3}
 \begin{gathered}
    \phi(\xi) = - \int_{B(\xi,2\delta)} \frac 1{8\pi^2} \chi(r)\, \log d_g(\xi,\eta)\, P_g\phi(\eta)\, dV_g(\eta) \\
    \qquad\qquad\qquad\qquad {} - \int_{B(\xi,2\delta)\setminus B(\xi,\delta)} E_1(\xi,\eta)\, \phi(\eta)\, e^{4w(\eta)}
    \, dV_g(\eta).
 \end{gathered}
\end{equation}
Written weakly, this reads:
\[
    P_g \left( -\frac 1{8\pi^2} \chi(r)\, \log d_g(\xi, \eta) \right) = \delta_{\xi} + \text{a bounded function},
\]
as claimed in Lemma~\ref{jan18l1}.

Finally we show $G(\xi,\eta)=G(\eta,\xi)$. This part of the proof
is similar to the proof in \cite[page 108]{aubin}. We include it
here for completeness.  We first claim:
\begin{equation}    \label{feb18e3}
    \int_M G(\eta,\xi) \, dV_g(\eta) = \text{constant}.
\end{equation}
Indeed, let $f(\xi)=\int_M G(\eta,\xi)\, dV_g(\eta)$, then for any
$\phi\in C^4(M)$,
\begin{multline*}
    \int_M f(\xi) P_{g,\xi} \phi(\xi)\, dV_g(\xi)
    = \int_M \left( \int_M G(\eta,\xi) P_{g,\xi} \phi(\xi)\, dV_g(\xi) \right)\, dV_g(\eta) \\
    = \int_M (\phi(\eta)-\bar \phi_g) \, dV_g(\eta) = 0.
\end{multline*}
Hence $P_gf=0$ in weakly, and since $\Ker(P_g)=\{\text{constants}\}$, we have
proved (\ref{feb18e3}). Next for any $\phi\in C^{4}(M)$,
\[
    \phi(\eta) - \bar\phi_g = \int_M G(\eta,\xi) P_{g,\xi} \phi(\xi)\, dV_g(\xi).
\]
hence for any $\psi\in C^4(M)$:
\begin{multline*}
    \int_M P_{g,\xi} \phi(\xi) \psi(\xi)\, dV_g(\xi) = \int_M \phi(\eta) P_{g,\eta} \psi(\eta) \, dV_g(\eta) \\
    = \int_M \bigl( \int_M G(\eta,\xi) P_{g,\xi} \phi(\xi)\, dV_g(\xi) + \bar\phi_g \bigr) P_{g,\eta} \psi(\eta)\, dV_g(\eta) \\
    = \int_M \bigl( \int_M G(\eta,\xi) P_{g,\eta} \psi(\eta)\, dV_g(\eta) \bigr) P_{g,\xi} \phi(\xi)\, dV_g(\xi).
\end{multline*}
Since $\Ker(P_g)=\{\text{constants}\}$, we obtain
\[
    \psi(\xi)=\int_MG(\eta,\xi)P_{g,\eta}\psi(\eta) \, dV_g(\eta) + \text{constant}.
\]
Thus, using also the definition of $G(\xi,\eta)$, we obtain
\[
    \int_M (G(\xi,\eta) - G(\eta,\xi)) P_{g,\eta} \psi(\eta) \, dV_g(\eta) = \text{constant}.
\]
Integrating with respect to $\xi$ and using~\eqref{intG=0} and~\eqref{feb18e3}, we see that the constant on the right hand side above is $0$.  Thus, again by~\eqref{P}
we have
\[
    G(\xi,\eta)-G(\eta,\xi) = C.
\]
where $C$ is independent of $\eta$.  Since the left-hand side is clearly antisymmetric with respect to $\xi$ and $\eta$, we obtain $G(\xi,\eta)=G(\eta,\xi)$. This completes the proof of Lemma \ref{jan18l1}.
\end{proof}

\section{Comparison between $d_{\gb}(y,z)$ and $|y-z|$}   \label{appdist}

In this subsection we establish the following estimate:
\begin{equation}    \label{jun27e10}
    D^j_y \bigl( \log |y-z| - \log d_{\breve g_k}(y,z) \bigr) = O(\eps_k^2|y|^{2-j}),
    \,\, \text{for $|z|<\frac{|y|}2$ and $|y|<\delta_1\eps_k^{-1}$,}
\end{equation}
for $j=1,2,3$. We recall that $\breve g_k=\eps_k^{-2}\phi_*g_k$ is the
blow-up metric, and we identify $y,z \in T_0 M$ with $\exp y$ and $\exp z$ respectively, where $\exp$ is the exponential map at the origin with respect to the metric $\gb_k$.  Thus, $d_{\gb_k}(y,z)=d_{\gb_k}(\exp y,\exp z)$.

Fix $y$ and consider
\[
    f(z) = \log |y-z|- \log d_{\gb_k}(y,z), \qquad \text{on $|z|<\frac 23|y|$},
\]
as a function of $z$.  We shall obtain the following estimate:
\begin{equation}    \label{jun27e11}
    |D^j_z f(z)| \le C\eps_k^2 |y|^{2-j},  \qquad \text{for $|z|<|y|/2$, $j=1,2,3.$}
\end{equation}
Then since $D_y^j f(z) = (-1)^j D_z^j f(z)$, \eqref{jun27e10}
follows.

Let $\breve R$, $\breve R_{ij}$, and $\breve R_{ijkl}$ respectively denote the
scalar, Ricci, and Riemann curvature of $\breve g$, where we suppress the subscript $k$.
By the definitions of $\breve g$ and $\breve R^i_{jlm}(z)$ one
obtains easily that
\[
 \nabla^j\breve R^i_{jml}(z)=O(\eps_k^{2+j}),\quad j=0,1,2.
\]
As a consequence, using $\breve R(0)=|\nabla \breve R(0)|=\breve
R_{ij}(0)=0$, we  have further:
\begin{equation}\label{jun27e2}
\breve R(z)=O(\eps_k^4|z|^2), \quad \breve R_{ij}(z)=
O(\eps_k^3|z|).
\end{equation}
We also note the following simple estimate on $d_{\gb}$:
\begin{equation}    \label{mar26e2}
    |D^j_z(\log d_{\gb}(y,z)) |\le C|y-z|^{-j}, \qquad j=1,2,3,4.
\end{equation}

We shall now derive an estimate on $\Delta^2_{\gb}f(z)$.  By~\eqref{jun27e1} and the definition of
$\gb$, we have:
\begin{equation}\label{jun27e3}
    P_{\gb,z} \log d_{\gb}(y,z) = O(\eps_k^4). \qquad |z|< \frac23|y|, \qquad |y| \le \delta_1\eps_k^{-1}.
\end{equation}
Next, using~\eqref{jun27e2} and~\eqref{mar26e2}, we can estimate the term:
\begin{multline*}
    (P_{\gb} - \Delta_{\gb}^2) \log d_{\gb}(y,z) \\
     = \partial_m \left( \gb^{mi} \left( \frac23
    \breve R(z) \gb_{ij} - 2 \breve R_{ij}(z) \right) \gb^{lj} \partial_j (\log d_{\gb}(y,z))\right)
    = O(\eps_k^3|y|^{-1}).
\end{multline*}
Combining this with~\eqref{jun27e3}, we get:
\begin{equation}    \label{jun27e5}
    \Delta_{\gb,z}^2 (\log d_{\gb}(y,z)) = O(\eps_k^3|y|^{-1}), \qquad |z| < \frac23|y|, \qquad |y|\le \delta_1\eps_k^{-1}.
 \end{equation}

Finally, we consider the term $\Delta_{\gb,z}^2(\log |y-z|)$.  Since $\Delta_z^2(\log |y-z|)=0$, it suffices to estimate $\Delta_{\gb,z}^2 - \Delta_z^2$.  For any function $u$, we have,
by direct computation:
\begin{multline}    \label{jun27e6}
    \Delta_{\gb,z}^2 u = \gb^{ab} \gb^{ij} \partial_{ijab} u + 2\partial_{ija} u (\partial_b \gb^{ab} \gb^{ij}
    + \gb^{ab} \partial_b \gb^{ij}) \\
    {} + \partial_{ij}u(\partial_a\gb^{ab}\partial_b\gb^{ij}+2\gb^{ai}\partial_{ab}\gb^{bj}+\gb^{ab}\partial_{ab}\gb^{ij}
    + \partial_a\gb^{ia}\partial_b\gb^{bj}) \\
    {} + \partial_j u (\partial_a \gb^{ab} \partial_{ib} \gb^{ij} + \gb^{ab} \partial_{iab} \gb^{ij}).
\end{multline}
where we used $\det(\gb)=1$.  Using the expansion of $\gb_{ij}(z)$:
\[
    \gb_{ij}(z) = \delta_{ij} + \frac13 \eps_k^2 \Rh_{pijq}(0) z^{ij} + O(\eps_k^3|z|^3).
\]
where $z^{ij}=z^iz^j$, and replacing $u$ by $\log |y-z|$ in~\eqref{jun27e6}, we obtain:
\[
    \Delta_{\gb,z}^2 (\log |y-z|) = O(\eps_k^2|y|^{-2}), \qquad |z|<\frac23|y|, \qquad |y|\le \delta_1\eps_k^{-1}
\]
and consequently
\begin{equation}    \label{jun27e8}
    \Delta_{\gb,z}^2 f(z) = O(\eps_k^2|y|^{-2}), \qquad |z|<\frac23|y|, \qquad |y|\le \delta_1\eps_k^{-1}.
\end{equation}

An estimate on the $L^{\infty}$ norm of $f(z)$ is easily obtained:
\begin{equation}    \label{feb13e3}
    d_{\gb}(y,z)=\int_z^y\sqrt{\gb_{ij}(t)x'_i(t)x'_j(t)}dt=|y-z|(1+O(\eps_k^2(|y|^2+|z|^2))).
\end{equation}
Note that here we didn't need the assumption: $|z|<\frac 23|y|$.  From~\eqref{feb13e3}, we clearly have
\begin{equation}    \label{jun27e9}
    f(z)=O(\eps_k^2|y|^2) \qquad |z|<\frac23 |y|, \qquad |y|<\delta_1\eps_k^{-1}.
 \end{equation}
With~\eqref{jun27e8} and~\eqref{jun27e9}, we can apply the standard rescaling argument and elliptic theory to derive~\eqref{jun27e11}.

\section{Conformal normal coordinates}  \label{appcoords}

In this subsection we list some well known facts for convenience.
Let $g$ be a metric on $M$ and let $p\in M$ be a point. We emphasize that $g$ and $u$ are unrelated to the corresponding quantities in other sections. Suppose that in normal coordinates around $p$ the metric $g$ satisfies $\det(g)=1$, i.e., those normal coordinates are conformal normal coordinates at $p$.  We will denote $p$ as $0$ and
consider the properties of $g$ in a neighborhood of $0$.

First, since $\det(g)=1$, the Laplacian is given by:
\begin{equation}    \label{nccDelta}
    \Delta_g u = \partial_j (g^{ij}\partial_j u) = \partial_j g^{ij} \partial_j u + g^{ij} \partial_{ij} u.
\end{equation}

The second term $P_g u - \Delta_g^2 u$ of the Paneitz operator is given by:
\begin{equation}    \label{nccPg}
    \div_g\left( \left( \frac23 R g - 2\Ric \right)\, du\, \right) =
    \partial_m\left( \left( \frac23 R g_{ij} - 2R_{ij} \right) \partial_l u\, g^{lj}\, g^{mi} \right).
\end{equation}

Further properties of conformal normal coordinates include:
\begin{gather}  \label{nov29e3}
    R_{ij}(0) = 0 \\[1ex]
    R_{ij,k}(0) + R_{jk,i}(0) + R_{ki,j}(0) = 0 \\[1ex]
    \nabla R(0) = 0 \\[1ex]
    \Delta R(0) = -\frac16 |W(0)|^2 \\
    R_{(ij,kl)}(0) + \frac29 R_{(pijm}(0) R_{pklm)}(0) = 0 \\[1ex]
    \label{nov29e2}
    R_{pijq,p}(0) = R_{iq,j}(0) - R_{ij,q}(0).
\end{gather}
Here, we use the round brackets as the customary notation for the symmetric part with respect to all indices not contracted within a certain scope.

In normal (not necessarily conformal normal) coordinates, there holds:
\begin{equation}    \label{oct22e0}
    g_{ab}(\xi) = \delta_{ab} + \frac13 R_{aijb}(0) \xi^{ij} + \frac16 R_{aijb,k}(0) \xi^{ijk} +
    O(r^4),\quad r=|\xi |
\end{equation}
and consequently, the same expansion holds for the inverse:
\begin{equation}    \label{oct22e1}
    g^{ab}(\xi) = \delta_{ab} - \frac13 R_{aijb}(0) \xi^{ij} - \frac16 R_{aijb,k}(0) \xi^{ijk} + O(r^4).
\end{equation}
Here $\xi^{ij}$ denotes $\xi^i\,\xi^j$, etc.

Taking a derivative:
\begin{equation}    \label{oct22e9}
    \partial_c g^{ab}(\xi) = -\frac23 R_{a(ci)b}(0) \xi^i -\frac16 \bigl( 2 R_{a(ci)b,j}(0) + R_{aijb,c}(0) \bigr)
     \xi^{ij} + O(r^3).
\end{equation}

Contracting over to $a$ and $c$ and using~\eqref{nov29e2}:
\begin{equation} \label{oct22e3}
    \partial_a g^{ab}(\xi) = -\frac16 \bigl( 2R_{ib,j}(0) - R_{ij,b}(0) \bigr) \xi^{ij} + O(r^3).
\end{equation}

Taking another derivative in~\eqref{oct22e9}:
\begin{multline}
    \partial_{cd} g^{ab}(\xi) = -\frac23 R_{a(cd)b}(0) \\
    -\frac13 \bigl(R_{a(cd)b,i}(0) + R_{iba(c,d)}(0) - R_{aib(c,d)}(0) \bigr) \xi^i + O(r^2).
\end{multline}

Contracting over $a$ and $c$ and using~\eqref{nov29e2} and~\eqref{nov29e3}:
\begin{equation} \label{nov29e1}
    \partial_{ad}g^{ab}(\xi) = \frac23R_{id,b}(0) \xi^i + O(r^2).
\end{equation}

\section{A Pohozaev Identity}   \label{appPohozaev}

In this appendix, we derive a Pohozaev identity for the equation
\begin{equation}    \label{Peq}
    P_g u + 2b = 2h e^{4u}.
\end{equation}
Throughout, we assume that $\det(g)=1$ over $\Omega$ which we take to be a ball centered at $0$.

First, multiplying the right hand side of~\eqref{Peq} by $\xi^i\partial_i u$ and integrating by parts, we have:
\[
    \int_{\Omega} 2he^{4u} \xi^i \partial_i u = \frac12 \int_{\partial\Omega} he^{4u} \xi^i \nu_i -
    \int_{\Omega}(2he^{4u} + \frac12\xi^i\partial_i he^{4u}),
\]
where $\nu_i$ is the unit normal to the boundary $\partial\Omega$.
Also note that we omitted to write here, as we will throughout
this appendix, the volume element $dV_g$ in the integral over
$\Omega$, and the area element $dA_g$ in the integral over
$\partial\Omega$.

Next, we consider the term $\Delta^2_g u$ on the right hand side of~\eqref{Peq}.  Again, multiplying by the same factor $\xi\cdot\nabla u$ and integrating by parts, we get:
\begin{multline*}
    \int_{\Omega} \Delta_g^2 u\, (\xi\cdot \nabla u)
    = \int_{\Omega} \partial_i \bigl( g^{ij} \partial_j(\Delta_g u) \bigr) \, \xi^k\partial_k u \\
    = \int_{\partial\Omega} g^{ij} \partial_j(\Delta_g u)\, \xi^k\partial_k u\, \nu_i
    - \int_{\Omega}g^{ij}\partial_j(\Delta_g u)\, \partial_i(\xi^k\partial_k u) \\
    = I_1 - \int_{\Omega} g^{ij} \partial_j (\Delta_g u)\, \partial_i u
    - \int_{\Omega} g^{ij} \partial_j (\Delta_g u)\, \xi^k \partial_{ik} u \\
    = I_1 - \int_{\partial\Omega} g^{ij} \Delta_g u\, \partial_i u\, \nu_j
    + \int_{\Omega} (\Delta_g u)^2
    - \int_{\partial\Omega} g^{ij} \Delta_g u\, \xi^k \partial_{ik} u\, \nu_j \\
    \qquad\qquad\qquad\qquad {} + \int_{\Omega} \Delta_g u\, \partial_j (g^{ij} \xi^k \partial_{ik} u) \\
    = I_1 - I_2 + \int_{\Omega} (\Delta_g u)^2 - I_3 + \int_{\Omega} g^{ij} \Delta_g u\, \partial_{ij} u
    + \int_{\Omega} \Delta_g u\, \xi^k \partial_j (g^{ij} \partial_{ik} u).
\end{multline*}
Here $I_1$, $I_2$ and $I_3$ are boundary integrals. In order to compute the last two terms we now note:
\begin{gather*}
    g^{ij} \partial_{ij} u = \Delta_g u - \partial_i g^{ij} \partial_j u \\
    \partial_j(g^{ij}\partial_{ik} u) = \partial_k (\Delta_g u) - \partial_{ik} g^{ij} \partial_j u - \partial_k g^{ij} \partial_{ij} u.
\end{gather*}
These imply:
\[
    \int_{\Omega} g^{ij} \Delta_g u\, \partial_{ij} u
    = \int_{\Omega}(\Delta_gu)^2 - \int_{\Omega} \Delta_g u\, \partial_i g^{ij} \partial_j u
    =\int_{\Omega} (\Delta_gu)^2 - B_3,
\]
and
\begin{multline*}
    \int_{\Omega} \Delta_g u \, \xi^k \partial_j (g^{ij}\partial_{ik} u) \\
    = \int_{\Omega} \xi^k \Delta_g u\, \partial_k(\Delta_g u)
    - \int_{\Omega} \xi^k \Delta_g u\, \partial_{ik} g^{ij} \partial_j u
    - \int_{\Omega} \xi^k \Delta_g u\, \partial_k g^{ij} \partial_{ij} u \\
    = \frac12 \int_{\Omega} \xi^k \partial_k\left( (\Delta_g u)^2 \right) - B_1 - B_2 \\
    = \frac12 \int_{\partial \Omega} (\xi\cdot\nu) (\Delta_g u)^2 - 2\int_{\Omega} (\Delta_g u)^2 - B_1 - B_2.
\end{multline*}
Thus, we conclude:
\[
    \int_{\Omega} \Delta_g^2 u\, (\xi\cdot \nabla u) = I_1 - I_2 - I_3 + \frac12 I_4 - B_1 - B_2 - B_3.
\]

Finally, we consider the second term in the Paneitz operator.  As before, we multiply by $\xi^k\partial_k u$ and integrate by parts:
\begin{multline*}
    \int_{\Omega} \partial_m \left( g^{mi} \left( \frac23 R(\xi)\, g_{ij}
        - 2R_{ij}(\xi) \right) g^{lj} \partial_l u  \right) \xi^k \partial_k u \\
    = \int_{\partial\Omega} g^{mi} \left( \frac23 R(\xi)\, g_{ij} - 2R_{ij}(\xi) \right) g^{lj} \partial_l u\, \xi^k \partial_k u\, \nu_m \\
    \qquad\qquad\qquad {}
    - \int_{\Omega} g^{mi} \left( \frac 23R(\xi)\,g_{ij}
        - 2R_{ij}(\xi) \right)\, g^{lj} \partial_l u\, (\partial_m u + \xi^k \partial_{mk} u) \\
    {} = C_1 - \int_{\Omega} g^{mi} \left( \frac23R(\xi)\,g_{ij} - 2R_{ij}(\xi) \right)
            g^{lj}\partial_l u\, (\partial_mu+\xi^k\partial_{mk}u)
\end{multline*}
Since $R_{ij}(0)=0$, $R(0)=0$, and $\partial_iR(0)=0$, we see that $C_1$ above can be written as
\[
    C_1=\int_{\partial\Omega} (-2R_{ij,l}(0)\, \partial_j u\, \partial_k u\,\xi^l\,\xi^k\,\nu_i +
    O(r^3)|Du|^2),\quad r=|\xi |.
\]
We also estimate the second term:
\begin{multline*}
    \int_{\Omega}g^{mi} \left( \frac23R(\xi)\,g_{ij} - 2R_{ij}(\xi) \right)\, \partial_l u \, g^{lj} (\partial_mu+\xi^k\partial_{mk}u) \\
    = \int_{\Omega} \bigl( (O(r^2) - 2R_{ij,s}(0)\, \xi^s)\, \partial_l u \, g^{lj} (\partial_i u + \xi^k\, \partial_{ik} u )
    + O(r^3)\,|Du|\,(|Du|+r|D^2u|) \bigr) \\
    = \int_{\Omega} \bigl( (O(r^2)-2R_{ij,s}(0)\, \xi^s)\, \partial_j u\, (\partial_i u + \xi^k\, \partial_{ik} u)
    + O(r^3)\,|Du|\,(|Du|+r|D^2u|) \bigr) \\
    = \int_{\Omega} \bigl( -2R_{ij,l}(0)\, \xi^l \partial_j u (\partial_i u + \xi^k\, \partial_{ik} u)
    + O(r^2) |Du|^2 + O(r^4)|D^2u| \bigr).
\end{multline*}

Finally, putting all the estimates together, we obtain the final form of our Pohozaev Identity:
\begin{multline}    \label{nov30e1}
    \int_{\Omega} \left( 2he^{4u} + \frac12\xi^i\,\partial_i h e^{4u} \right)\\
    = \int_{\partial\Omega} \left( \frac12 h e^u \xi^i\, \nu_i - g^{ij} \partial_i (\Delta_gu)\, \partial_k u \, \xi^k\, \nu_j
    + g^{ij} \, \Delta_g u\, \partial_i u\, \nu_j \right. \\
    \left. \qquad\qquad\qquad\qquad\qquad {} + g^{ij}\,\Delta_g u\, \xi^k\, \partial_{ik} u\, \nu_j
    - \frac12\,(\Delta_gu)^2\, \xi^i\, \nu_i\right) \\
    + \int_{\Omega} \left( \Delta_g u\, \partial_i\, g^{ij}\, \partial_j u + \xi^k\, \Delta_g u\, \partial_{ik}\, g^{ij}\, \partial_j u
    + \xi^k\, \Delta_g u\, \partial_k g^{ij}\, \partial_{ij} u - 2 b\, \xi^i\, \partial_i u \right) \\
    {} + 2 \int_{\partial \Omega} \left( R_{ij,l}(0)\, \partial_j u\, \partial_k u\, \xi^l\, \xi^k\, \nu_i
    + O(r^3)|Du|^2 \right) \\
    - \int_{\Omega} \bigl( 2 R_{ij,l}(0) ( \partial_j u\, \partial_i u\, \xi^l
    + \partial_j u\, \partial_{ik} u\, \xi^k\, \xi^l) + O(r^2)|Du|^2+O(r^4)|D^2u| \bigr).
\end{multline}


\begin{thebibliography}{99}

\newcommand{\aut}[1]{{\sc #1},}
\newcommand{\tit}[1]{{\em #1\/},}
\newcommand{\vol}[1]{{\bf #1}}
\newcommand{\yr}[1]{(#1)}
\newcommand{\pp}[2]{#1--#2}

\bibitem{adimurthi}
    \aut{Adimurthi, F.~Robert and M.~Struwe}
    \tit{Concentration phenomena for Liouville's equation in dimension four}
    J.\ Eur.\ Math.\ Soc.\
    \vol{8}
    \yr{2006} no.~2,
    \pp{171}{180}.

\bibitem{aubiny}
    \aut{T.~Aubin}
    \tit{\'Equations diff\'erentielles non lin\'eaires et probl\`eme de Yamabe concernant la courbure scalaire}
    J.\ Math.\ Pures Appl.\
    \vol{55}
    \yr{1976},
    \pp{269}{296}.

\bibitem{aubin}
    \aut{T.~Aubin}
    \tit{Some nonlinear problems in Riemannian geometry}
    Springer Monographs in Mathematics, Springer-Verlag, Berlin, 1998.

\bibitem{brendle}
    \aut{S.~Brendle}
    \tit{Convergence of the Q-curvature flow on S4}
    Adv.\ Math.\
    \vol{205}
    \yr{2006} no.~1,
    \pp{1}{32}.

\bibitem{browder}
    \aut{F.~E.~Browder}
    \tit{On the regularity properties of solutions of elliptic differential equations}
    Comm.\ Pure Appl.\ Math.\
    \vol{9}
    \yr{1956},
    \pp{351}{361}.

\bibitem{cao}
    \aut{J.~G.~Cao}
    \tit{The existence of generalized isothermal coordinates for higher-dimensional Riemannian manifolds}
    Trans.\ Amer.\ Math.\ Soc.\
    \vol{324}
    \yr{1991} no.~2,
    \pp{901}{920}.

\bibitem{chang1}
    \aut{S.~A.~Chang}
    \tit{On a fourth-order partial differential equation in conformal geometry}
    in \tit{Harmonic analysis and partial differential equations (Chicago, IL, 1996)} 127--150,
    Chicago Lectures in Math., Univ.\ Chicago Press, Chicago, IL, 1999.

\bibitem{changqingyang}
    \aut{S.~A.~Chang, J.~Qing and P.~Yang}
    \tit{Some progress in conformal geometry}
    SIGMA Symmetry Integrability Geom.\ Methods Appl.\
    \vol{3}
    \yr{2007},
    Paper 122.

\bibitem{changyang1}
    \aut{S.~A.~Chang and P.~Yang}
    \tit{Extremal metrics of zeta function determinants on $4$-manifolds}
    Annals of Math
    \vol{142}
    \yr{1995},
    \pp{171}{212}.

\bibitem{changyang2}
    \aut{S.~A.~Chang and P.~Yang}
    \tit{On a fourth order curvature invariant}
    in \tit{Spectral problems in geometry and arithmetic (Iowa City, IA, 1997)} 9--28,
    Contemp.\ Math.\ \vol{237}, Amer.\ Math.\ Soc., Providence, RI, 1999.

\bibitem{christodoulou}
    \aut{D.~Christodoulou and S.~Klainerman}
    \tit{The Global Nonlinear Stabililty of Minkowski Space}
    Princeton Mathematical Series, No.~41, Princeton University Press, Princeton, 1993.

\bibitem{djadli1}
    \aut{Z.~Djadli and A.~Malchiodi}
    \tit{A fourth order uniformization theorem on some four manifolds with large total $Q$-curvature}
    C.\ R.\ Math.\ Acad.\ Sci.\ Paris
    \vol{340}
    \yr{2005} no.~5,
    \pp{341}{346}.

\bibitem{djadli2}
    \aut{Z.~Djadli and A.~Malchiodi}
    \tit{Existence of conformal metrics with constant Q-curvature}
    Ann.\ of Math., to appear.

\bibitem{nirenberg2}
    \aut{A.~Douglis and L.~Nirenber}
    \tit{Interior estimates for elliptic systems of partial differential equations}
    Comm.\ Pure Appl.\ Math.\
    \vol{8}
    \yr{1955}
    \pp{503}{538}.

\bibitem{druet1}
    \aut{O.~Druet and F.~Robert}
    \tit{Bubbling phenomena for fourth-order four-dimensional PDEs with exponential growth}
    Proc.\ Amer.\ Math.\ Soc.\
    \vol{134}
    \yr{2006} no.~3,
    \pp{897}{908}.

\bibitem{fefferman}
    \aut{C.~Fefferman and K.~Hirachi}
    \tit{Ambient metric construction of $Q$-curvature in conformal and CR-geometries}
    Math.\ Res.\ Lett.\
    \vol{10}
    \yr{2003} no.~5-6,
    \pp{819}{831}.

\bibitem{gursky}
    \aut{M.~Gursky}
    \tit{The principle eigenvalue of a conformally invariant differential operator, with an application to semi-linear elliptic PDE}
    Comm.\ Math.\ Phys.\
    \vol{207}
    \yr{1999},
    \pp{131}{143}.

\bibitem{gurskyviaclovsky}
    \aut{M.~Gursky and J.~Viaclovsky}
    \tit{A fully nonlinear equation on four-manifolds with positive scalar curvature}
    J.\ Differential Geom.
    \vol{63}
    \yr{2003} no. 1,
    \pp{131}{154}.

\bibitem{gurskyviaclovsky2}
    \aut{M.~Gursky and J.~Viaclovsky}
    \title{A fully nonlinear equation on four-manifolds with positive scalar curvature}
    J.\ Differential Geom.\
    \vol{63}
    \yr{2003} no.~1,
    \pp{131}{154}.


\bibitem{koselev}
    \aut{A.~I.~Ko\u selev}
    \tit{On the boundedness in $L\sb{p}$ of derivatives of solutions of elliptic equations and elliptic systems} (Russian) Dokl.\ Akad.\ Nauk SSSR (N.S.)
    \vol{116}
    \yr{1957}
    \yr{542}{544}.

\bibitem{leeparker}
    \aut{J.~M.~Lee and T.~H.~Parker}
    \tit{The Yamabe problem}
    Bull.\ Amer.\ Math.\ Soc.\ (N.S.)
    \vol{17}
    \yr{1987} no.~1,
    \pp{37}{91}.

\bibitem{lililiu}
    \aut{J.~Li, Y.~Li and P.~Liu}
    \tit{The $Q$-curvature on a $4$-dimensional Riemannian manifold $(M,g)$ with $\int_M Q\,dV_g = 8\pi^2$}
    arXiv, math/0608543.

\bibitem{lin1}
    \aut{C.~S.~Lin}
    \tit{A classification of solutions of a conformally invariant fourth order equation in $\R^n$}
    Comment.\ Math.\ Helv.\
    \vol{73}
    \yr{1998} no.~2,
    \pp{206}{231}.

\bibitem{linwei1}
    \aut{C.~S.~Lin and J.~C.~Wei}
    \tit{Sharp Estimates For Bubbling Solutions of A Fourth Order Mean Field Equation}
     Ann.\ Sc.\ Norm.\ Super.\ Pisa Cl.\ Sci., to appear.

\bibitem{mal1}
    \aut{A.~Malchiodi}
    \tit{Compactness of solutions to some geometric fourth-order equations}
    J.\ Reine Angew.\ Math.\
    \vol{594}
    \yr{2006},
    \pp{137}{174}.

\bibitem{malchiodi}
    \aut{A.~Malchiodi}
    \tit{Conformal metrics with constant $Q$-curvature}
    SIGMA Symmetry Integrability Geom.\ Methods Appl.\
    \vol{3}
    \yr{2007},
    Paper 120.

\bibitem{malchiodi2}
    \aut{A.~Malchiodi and M.~Struwe}
    \tit{$Q$-curvature flow on $S^4$}
    J.\ Differential Geom.\
    \vol{73}
    \vol{2006} no.~1,
    \pp{1}{44}.

\bibitem{ndiaye}
    \aut{C.~B.~Ndiaye}
    \tit{Constant Q-curvature metrics in arbitrary dimension}
    J.\ Funct.\ Anal.\
    \vol{251}
    \yr{2007} no.~1,
    \pp{1}{58}.

\bibitem{nirenberg1}
    \aut{L. Nirenberg}
    \tit{Estimates and existence of solutions of elliptic equations}
    Comm.\ Pure Appl.\ Math.\
    \vol{9}
    \yr{1956},
    \pp{509}{529}.

\bibitem{qing}
    \aut{J.~Qing and D.~Raske}
    \tit{Compactness for conformal metrics with constant $Q$ curvature on locally conformally flat manifolds}
    Calc.\ Var.\ Partial Differential Equations
    \vol{26}
    \yr{2006} no.~3,
    \pp{343}{356}.

\bibitem{rwei}
    \aut{F.~Robert and J.~C.~Wei}
    \tit{Asymptotic behavior of a fourth order mean field equation with Dirichlet boundary condition}
    Indiana Univ.\ Math.\ J., to appear.

\bibitem{schoen}
    \aut{R.~Schoen}
    \tit{Conformal deformation of a Riemannian metric to constant scalar curvature}
    J.\ Differential Geom.\
    \vol{20}
    \yr{1984},
    \pp{479}{495}.

\bibitem{struwe}
    \aut{M.~Struwe}
    \tit{Quantization for a fourth order equation with critical exponential growth}
    (English summary)
    Math.\ Z.\
    \vol{256} no.~2
    \yr{2007}
    \pp{397}{424}.


\bibitem{trudinger}
    \aut{N. Trudinger}
    \tit{Remarks concerning the conformal deformation of Riemannian structures on compact manifolds}
    Ann.\ Scuola Norm.\ Sup.\ Pisa Cl.\ Sci.\ (3)
    \vol{22}
    \yr{1968},
    \pp{265}{274}.

\bibitem{wei}
    \aut{J.~Wei}
    \tit{Asymptotic behavior of a nonlinear fourth order eigenvalue problem}
    Comm.\ Partial Differential Equations
    \vol{21}
    \yr{1996} no.~9-10,
    \pp{1451}{1467}.

\bibitem{wei2}
    \aut{J.~Wei and X.~Xu}
    \tit{On conformal deformations of metrics on $S^n$}
    J.\ Funct.\ Anal.\
    \vol{157}
    \yr{1998} no.~1,
    \pp{292}{325.}.

\bibitem{weinstein}
    \aut{G.~Weinstein}
    \tit{The Poincar\'e Uniformization Theorem}
    unpublished note,
    \url{http://www.math.uab.edu/weinstei/notes/poincare.pdf}.

\bibitem{xu}
    \aut{Y.~Xu}
    \tit{An estimate on the blowing-up solutions of a fourth-order equation}
    J.\ Funct.\ Anal.\
    \vol{251}
    \yr{2007} no.~1,
    \pp{360}{375}.

\bibitem{yamabe}
    \aut{H.~Yamabe}
    \tit{On a deformation of Riemannian structures on compact manifolds}
    Osaka Math.\ J.\
    \vol{12}
    \yr{1960}
    \pp{21}{37}.

\end{thebibliography}
\end{document}